\documentclass[12pt]{article}

\usepackage{mathptmx,amsmath,amssymb,amsfonts,amsthm}
\usepackage{graphicx,color,xcolor,pictex,epstopdf,url,algorithm,algorithmic,caption,endnotes}

\usepackage{etoolbox} % another kludge to over-ride Latex defaults
\usepackage{wrapfig} %
\usepackage{mathtools}
\makeatletter
\patchcmd{\@makechapterhead}{\large}{\normalsize}{}{}% for \chapter
\patchcmd{\@makeschapterhead}{\normalsize}{\normalsize}{}{}% for \chapter*
\makeatother
\usepackage[textwidth=6.0in,textheight=8.9in,left=0.75in,right=0.75in,top=0.75in,bottom=0.75in]{geometry}
\usepackage[section]{placeins}
\makeatletter
\g@addto@macro\normalsize{\setlength\abovedisplayskip{4pt}}
\g@addto@macro\normalsize{\setlength\belowdisplayskip{4pt}}
\makeatother
\newtheorem{theorem}{Theorem}[section]

\newtheorem{corollary}{Corollary}[section]
\newtheorem{proposition}{Proposition}[section]

\newtheorem{remark}{Remark}[section]
\usepackage{parskip}
\let\oldref\ref
\renewcommand{\ref}[1]{(\oldref{#1})}  % stupid kludge for laTeX
\renewcommand{\eqref}[1]{(\oldref{#1})} %BB03dec18 

\newbox\boxaddrone \newbox\boxaddrtwo

\input colordvi
\font\tenrm=cmr10
\font\teni=cmmi10 \skewchar\teni='177
\font\tensy=cmsy10 \skewchar\tensy='60
\font\tenex=cmex10
\font\tenit=cmti10
\font\tensl=cmsl10
\font\tenbf=cmbx10
\font\tentt=cmtt10
\font\ninerm=cmr9
\font\ninei=cmmi9 \skewchar\ninei='177
\font\ninesy=cmsy9 \skewchar\ninesy='60
\font\nineit=cmti9
\font\ninesl=cmsl9
\font\ninebf=cmbx9
\font\ninett=cmtt9
\font\eightrm=cmr8
\font\eighti=cmmi8 \skewchar\eighti='177
\font\eightsy=cmsy8 \skewchar\eightsy='60
\font\eightit=cmti8
\font\eightsl=cmsl8
\font\eightbf=cmbx8
\font\eighttt=cmtt8
\font\sevenrm=cmr7
\font\seveni=cmmi7 \skewchar\seveni='177
\font\sevensy=cmsy7 \skewchar\sevensy='60
\font\sevenbf=cmbx7
\font\sevenit=cmmi7
\font\sevensl=cmmi7
\font\seventt=cmr7
\font\sixrm=cmr6
\font\sixi=cmmi6 \skewchar\sixi='177
\font\sixsy=cmsy6 \skewchar\sixsy='60
\font\sixbf=cmbx6
\font\fiverm=cmr5
\font\fivei=cmmi5 \skewchar\fivei='177
\font\fivesy=cmsy5 \skewchar\fivesy='60
\font\fivebf=cmbx5
\def\tenpoint{\def\rm{\fam0\tenrm}%
        \textfont0=\tenrm \scriptfont0=\sevenrm \scriptscriptfont0=\fiverm
        \textfont1=\teni \scriptfont1=\seveni \scriptscriptfont1=\fivei
        \textfont2=\tensy \scriptfont2=\sevensy \scriptscriptfont2=\fivesy
        \textfont3=\tenex \scriptfont3=\tenex \scriptscriptfont3=\tenex
        \def\it{\fam\itfam\tenit}%
        \textfont\itfam=\tenit
        \def\sl{\fam\slfam\tensl}%
        \textfont\slfam=\tensl
        \def\bf{\fam\bffam\tenbf}%
        \textfont\bffam=\tenbf \scriptfont\bffam=\sevenbf
                \scriptscriptfont\bffam=\fivebf
        \def\tt{\fam\ttfam\tentt}%
        \textfont\ttfam=\tentt
        \normalbaselineskip=12pt%
        \let\sc=\eightrm        % Small caps
        \setbox\strutbox=\hbox{\vrule height8.5pt depth3.5pt width0pt}%
        \normalbaselines\rm}
\def\ninepoint{\def\rm{\fam0\ninerm}%
        \textfont0=\ninerm \scriptfont0=\sixrm \scriptscriptfont0=\fiverm
        \textfont1=\ninei \scriptfont1=\sixi \scriptscriptfont1=\fivei
        \textfont2=\ninesy \scriptfont2=\sixsy \scriptscriptfont2=\fivesy
        \textfont3=\tenex \scriptfont3=\tenex \scriptscriptfont3=\tenex
        \def\it{\fam\itfam\nineit}%
        \textfont\itfam=\nineit
        \def\sl{\fam\slfam\ninesl}%
        \textfont\slfam=\ninesl
        \def\bf{\fam\bffam\ninebf}%
        \textfont\bffam=\ninebf \scriptfont\bffam=\sixbf
                \scriptscriptfont\bffam=\fivebf
        \def\tt{\fam\ttfam\ninett}%
        \textfont\ttfam=\ninett
        \normalbaselineskip=11pt%
        \let\sc=\sevenrm        % Small caps
        \setbox\strutbox=\hbox{\vrule height8pt depth3pt width0pt}%
        \normalbaselines\rm}
\def\eightpoint{\def\rm{\fam0\eightrm}%
        \textfont0=\eightrm \scriptfont0=\sixrm \scriptscriptfont0=\fiverm
        \textfont1=\eighti \scriptfont1=\sixi \scriptscriptfont1=\fivei
        \textfont2=\eightsy \scriptfont2=\sixsy \scriptscriptfont2=\fivesy
        \textfont3=\tenex \scriptfont3=\tenex \scriptscriptfont3=\tenex
        \def\it{\fam\itfam\eightit}%
        \textfont\itfam=\eightit
        \def\sl{\fam\slfam\eightsl}%
        \textfont\slfam=\eightsl
        \def\bf{\fam\bffam\eightbf}%
        \textfont\bffam=\eightbf \scriptfont\bffam=\sixbf
                \scriptscriptfont\bffam=\fivebf
        \def\tt{\fam\ttfam\eighttt}%
        \textfont\ttfam=\eighttt
        \normalbaselineskip=9pt%
        \let\sc=\sixrm  % Small caps
        \setbox\strutbox=\hbox{\vrule height7pt depth2pt width0pt}%
        \normalbaselines\rm}
\def\sevenpoint{\def\rm{\fam0\sevenrm}%
        \textfont0=\sevenrm \scriptfont0=\fiverm \scriptscriptfont0=\fiverm
        \textfont1=\seveni \scriptfont1=\fivei \scriptscriptfont1=\fivei
        \textfont2=\sevensy \scriptfont2=\fivesy \scriptscriptfont2=\fivesy
        \textfont3=\tenex \scriptfont3=\tenex \scriptscriptfont3=\tenex
        \def\it{\fam\itfam\sevenit}%
        \textfont\itfam=\sevenit
        \def\sl{\fam\slfam\sevensl}%
        \textfont\slfam=\sevensl
        \def\bf{\fam\bffam\sevenbf}%
        \textfont\bffam=\sevenbf \scriptfont\bffam=\fivebf
                \scriptscriptfont\bffam=\fivebf
        \def\tt{\fam\ttfam\seventt}%
        \textfont\ttfam=\seventt
        \normalbaselineskip=8pt%
        \let\sc=\fiverm  % Small caps
        \setbox\strutbox=\hbox{\vrule height6pt depth2pt width0pt}%
        \normalbaselines\rm}

\input pictex
\newdimen\xfiglen \newdimen\yfiglen  \newdimen\textlen

\def\ul#1{\underline{#1}}
\def\ol#1{\overline{#1}}

\def\R{\mathbb{R}}
\def\N{\mathbb{N}}

\begin{document}

\title{
On uniqueness and reconstruction of a nonlinear diffusion term in a parabolic equation}
\author{Barbara Kaltenbacher\footnote{
Department of Mathematics,
Alpen-Adria-Universit\"at Klagenfurt.
barbara.kaltenbacher@aau.at.}
\and
William Rundell\footnote{
Department of Mathematics,
Texas A\&M University,
%College Station,
Texas 77843. % USA.
rundell@math.tamu.edu}
}
\date{\vskip-3ex}
 \maketitle  

\begin{abstract}
The problem of recovering coefficients in a diffusion equation is one of the
basic inverse problems. Perhaps the most important term is
the one that couples the length and time scales and is often
referred to as {\it the\/} diffusion coefficient $a$ in 
$u_t - \nabla(a\nabla u) = f$. In this paper we seek the unknown $a$ assuming
that $a=a(u)$ depends only on the value of the solution at a given point.
Such diffusion models are the basic of a wide range of physical phenomena
such as nonlinear heat conduction, chemical mixing and population dynamics.
We shall look at two types of overposed data in order to effect recovery of
$a(u)$: the value of a time trace $u(x_0,t)$ for some fixed point $x_0$ on
the boundary of the region $\Omega$; or the value of $u$ on an interior
curve $\Sigma$ lying within $\Omega$. As examples, these might represent a 
temperature measurement on the boundary or a census of the population
in some subset of $\Omega$ taken at a fixed time $T>0$.
In the latter case we shall show a uniqueness result that leads to a
constructive method for recovery of $a$.
Indeed, for both types of measured data we shall show reconstructions based
on the iterative algorithms developed in the paper.
\end{abstract}

\leftline{\small \qquad\quad
}
\smallskip
\leftline{\small \qquad\quad 
\textbf{{\textsc ams}} {\bf classification:} 35R30, 35K15, 35K58, 80A23.}

%%%%%%%%%%%%%%%%%%%%%%%%%%%%%%%%%%%%%%%%%%%
\section{Introduction}

The setting is in a bounded, simply connected domain
$\Omega\subset \mathbb{R}^d$
with smooth ($C^2$) boundary $\partial\Omega$ and the problem is to recover
the conductivity coefficient  $a(u)$ in the reaction diffusion equation
\begin{eqnarray}
&&u_t-\nabla(a(u)\nabla u)= r(x,t,u) \quad t\in(0,T)\,, \quad u(0)=u_0\label{eqn:u-nl}
\end{eqnarray}
subject to
(nonlinear) impedance or Dirichlet ($\gamma=\infty$) boundary conditions 
\begin{equation}\label{eqn:bndy}
a(u)\partial_\nu u+\gamma u = b(x,t), \mbox{ or } u = d(x,t) \quad x\in\partial\Omega,\quad t\in(0,T)
\end{equation}
where we assume that the forcing term $r(x,t,u)$ is known.

Since for a given $a(u)$ within a suitable class,
the pair \eqref{eqn:u-nl} and \eqref{eqn:bndy}
allows a unique determination of the solution $u$,
we have to give additional (overposed) data in order to recover $a$.
This will take the form of either:
observations along a curve $\omega\subset\Omega$ for some fixed time $T$
(final time measurements)
\begin{equation}\label{eqn:fiti}
g(x)=u(x,T), \quad x\in \omega\subseteq\Omega.
\end{equation}
or 
%BK 2020-10-16 from 
time trace observations 
\begin{equation}\label{eqn:titr}
h(t)=u(x_0,t), \quad t\in (0,T)
\end{equation}
for some $x_0\in\overline{\Omega}\subseteq\R^d$.
In the latter case, typically $x_0$ will be a boundary point on
$\partial\Omega$.

There is a parallel equation to \eqref{eqn:u-nl} which in its simplest form is
\begin{equation}\label{eqn:phiu}
u_t- \triangle \phi(u) = 0 \quad t\in(0,T)\,, 
\end{equation}
and thus looking at only the principal part of the operator gives
the operator in \eqref{eqn:u-nl} with 
%BK 2020-10-16 $a(u) = \nabla \phi(u)$.
$a(u) = \phi'(u)$.

In many applications this is the preferred form, perhaps the best known
case is when $\phi(u) = u^m$ obtaining the {\it Porous Medium Equation,}
\cite{Vazquez:2007}.

The heat equation can be derived from a basic random walk model
whereby at fixed times $dt$ jumps of length $dx$ are made in a random direction;
that is simple Brownian motion. The coefficient $a$ appears as the coupling
constant between these length-time scales and may be a function of the position
$x$ -- or a coefficient $a(u)$ meaning that in the random walk model the jump
length depends on the value or density at the current location.
It plays the role of a nonlinear thermal conductivity in
heat conduction or of a nonlinear diffusion in wide variety of areas
including applications to population dynamics. 

The conductivity of materials depends on temperature and sometimes quite
strongly so that the simple assumption that it be constant is only
valid over an often fairly narrow range.
The thermal conductivity for pure  metals depends on the motion
of free electrons.
The molecular vibrations increase with temperature thus decreasing the mean
free path of molecules and in turn obstruct the flow of free electrons,
resulting in a  reduction of conductivity.
However, for most non-metals the opposite effect occurs and conductivity
often rises with temperature.
Alloys and composite materials can have more complex behaviour and
for most materials at extreme temperatures the graph of $a(u)$
may certainly not be monotonic.

Both the  forms \eqref{eqn:u-nl} and \eqref{eqn:phiu}
occur frequently in the life sciences and  many from ecology
are based on the seminal paper of Skellam, \cite{Skellam:1951}.
For example in \cite{Turchin:1989}, modelling aggregative movement in
populations where $\phi(u)$ was a cubic and the associated $a(u)$
was a constant plus a logistical growth term.
In \cite{GurtinMacCamy:1977,AudritoVazquez:2017} the doubly nonlinear
model with $r(x,t,u) = f(u)$ was modelling instead dispersive behaviour
of populations.
In \cite{GurtinMacCamy:1977}, a similar nonlinear
model was used where the reaction term $f(u)$ is based on the Fisher model of
the logistic quadratic nonlinearity.

Note that we assume that $r(x,t,u)$ is known although in the case of
%BK 2020-10-16 $a(u)=a(x)$
$a=a(x)$ 
and unknown $f(u)$, recovery of $f$ from over-posed data similar
to that considered here was studied in~\cite{KaltenbacherRundell:2019c}.

The presence of $a$ in the boundary condition~\eqref{eqn:bndy} 
follows from the fact we have imposed a condition on the thermal
flux rather than on the gradient of $u$ itself.
We will use this fact in a key way in our reconstruction algorithm
for time trace data where the representation will lead naturally to a
fixed point equation for $a(u)$.
This latter problem was studied in one space dimension by
Cannon and DuChateau~\cite{CannonDuChateau:1980}
and DuChateau~\cite{DuChateau:1981}
using the transformation 
%BK 2020-10-16 $b=T(a)$ 
$A(s)=\int_0^s a(r)\, dr$, $b(s)=(a(A^{-1}(s)))^{-1}$, $v(x,t)=A(u(x,t))$,
taking $u_t - \bigl(a(u)u_x\bigr)_x = 0$
into $b(u)u_t - u_{xx} = 0$ and then showing uniqueness using
monotonicity arguments.

In the next section we will consider fixed point schemes for the
iterative reconstruction of $a$ in these two cases of overposed data.
%BK 2020-10-16: eqref -> oldref 4 times 
In section~\oldref{sect:forwards} we shall prove some results on the forwards
problem for \eqref{eqn:u-nl} and \eqref{eqn:bndy} that will be needed
in the analysis of the inverse problems.
In section~\oldref{sec:fiti_contractivity} we shall show that the map developed in
section~\oldref{sect:setup} for final time data is contractive in a suitable
setting and leads directly to a uniqueness result.
The final section shows some reconstructions based on the above analysis.

\section{The problem set up}\label{sect:setup}

We will consider fixed point schemes for the iterative reconstruction of $a$ in these two cases of overposed data.
As will be seen, the recovery process both numerically and analytically is
quite different.
In some ways this is to be expected for although space and time is
intrinsically connected in a parabolic equation, the nonlinearity
has its own effects that can play out in a different manner within
the interior of $\Omega$ and on its boundary $\partial\Omega$.
In both cases we must arrange the prescribed data to ensure that
the range of $u$ over the measurement region contains the entire range of
$u$ over $\Omega\times (0,T)$.
This is usually unnecessary in the case of spatially dependent coefficients
and is one significant factor which sets recovery within equations such as
\eqref{eqn:u-nl} apart.
What this means is the ideal situation is the ability to set up an experiment
allowing adherence to these conditions and using known properties of parabolic
equations (such as the maximum principle) to show the viability.
This of course limits the ability to recover $a(u)$ in, say ecological models,
where one has usually to measure just what is found.

\subsection{Final time data}\label{subsec:fiti}

In \eqref{eqn:fiti}, the observation domain $\omega$ is supposed to be chosen
such that there exists a curve $\Sigma\subseteq\mbox{int}(\omega)$,
(so that by differentiation of the data we also know $\nabla g$,
$\triangle g$ on $\Sigma$) that can be parametrized in such a way that
\begin{equation}\label{eqn:Sigma}
\Sigma =\{\vec{x}(\sigma)\, : \sigma\in [0,1]\}\,, \quad 
|(g\circ\vec{x})'(\sigma)| = |\nabla g(\vec{x}(\sigma))\cdot \vec{x}'(\sigma)|\geq\kappa>0\,.
\end{equation}

This allows us to define the inverse operator 
\begin{equation}\label{eqn:phi}
\phi:= (g\circ\vec{x})^{-1}:g(\Sigma)\to[0,1].
\end{equation}
In the 1-d case $\Omega=\Sigma=(0,1)$, \eqref{eqn:Sigma} simplifies to the strict monotonicity assumption $|g_x(x)|\geq\kappa>0$.

Moreover we assume the range of $g$ to contain the range of $u_{act}$
\begin{equation}\label{eqn:range_fiti}
J:=[\ul{u},\ol{u}]=g(\Sigma)\supseteq u_{act}(\Omega\times(0,T)) 
\end{equation}
for an exact solution $(a_{act},u_{act})$ of the inverse problem 
\eqref{eqn:u-nl}, \eqref{eqn:bndy}, \eqref{eqn:fiti}.
Similarly to 
%BK 2020-10-16 BB5, 
\cite{KaltenbacherRundell:2020a}
we will then enforce $J$ as the domain of the reconstructions, which are set to constant values outside $J$.

Projecting the {\sc pde}  in \eqref{eqn:u-nl} on the observation manifold $\Sigma\times\{T\}$ yields the fixed point iteration
\begin{equation}\label{eqn:proj}
\nabla(a_{k+1}(g)\nabla g) = D_t u(\cdot,T;a_k)-r(T) \mbox{ on }\Sigma
\end{equation}
and further inverting the curve parametrization, cf. \eqref{eqn:Sigma}, we get
\begin{equation}\label{eqn:ODEa+}
\begin{aligned}
&&\hspace*{-2cm}
a_{k+1}'(\tau)|\nabla g(\vec{x}(\phi(\tau)))|^2 + a_{k+1}(\tau)\triangle g(\vec{x}(\phi(\tau))) = D_t u(\vec{x}(\phi(\tau)),T;a_k)-r(\vec{x}(\phi(\tau)),T)  \\ 
&&\mbox{ for all }\tau\in[\ul{u},\ol{u}].
\end{aligned}
\end{equation}
This {\sc ode}, together with a given value of $a$ at some fixed point
allows us to uniquely determine the update $a_{k+1}$
%BK 2020-10-16  
inside $J$. We set $a_{k+1}(\tau)=a_{le}(\tau)$ for $\tau<\ul{u}$, $a_{k+1}(\tau)=a_{ri}(\tau)$ for $\tau>\ol{u}$, where $a_{le}$, $a_{ri}$ are the -- assumed known -- values of $a_{act}$ outside $J$.
More precisely, to preserve $C^{2,\alpha}$ smoothness of the iterates over all of $\mathbb{R}$, we use the metric projection of the iterates onto $\{a\in H^s(\mathbb{R})\, : \, a(\tau)=a_{le}(\tau)\mbox{ for }\tau<\ul{u}, \ a(\tau)=a_{ri}(\tau)\mbox{ for }\tau>\ol{u}\}$
with $s>\frac{d}{2}+2+\alpha$ so that $H^s(\mathbb{R})$ continuously embeds into $C^{2,\alpha}(\mathbb{R})$, cf. \cite{KaltenbacherRundell:2019c,KaltenbacherRundell:2020b}.

In the 1-d case $\Omega=\Sigma=(0,1)$ this can simply be solved by integrating
\eqref{eqn:proj} from $0$ to $x$ and dividing by $g_x$, then replacing
$-a_{k+1}(g(0))g_x(0)$ by $b(0,T)-\gamma g(0)$ (cf. \eqref{eqn:bndy})
\begin{equation}\label{eqn:a+1-d}
a_{k+1}(g(x))=\frac{\int_0^x (D_t u(\xi,T;a_k)-r(\xi,T))\, d\xi - b(0,T) + \gamma g(0)}{g_x(x)}
\end{equation}
 In higher space dimensions, assuming that $x_0=\vec{x}(0)\in\Sigma\cap\partial\Omega$ and $\partial_\nu g(x_0)\not=0$ so that we can replace $a_{k+1}(0)$ by $a_{act}(g(x_0))=\frac{b(x_0,T)-\gamma g(x_0)}{\partial_\nu g(x_0)} =: a^0$ (cf. \eqref{eqn:bndy}), we get from \eqref{eqn:ODEa+}
\begin{equation}\label{eqn:a+h-d}
\begin{aligned}
a_{k+1}(\tau)=&\exp\left(-\int_0^\tau\frac{\triangle g(\vec{x}(\phi(\sigma)))}{|\nabla g(\vec{x}(\phi(\sigma)))|^2}\, d\sigma\right)
a^0\\
&+ \int_0^\tau \frac{D_t u(\vec{x}(\phi(\sigma)),T;a_k)-r(\vec{x}(\phi(\sigma)),T)}{|\nabla g(\vec{x}(\phi(\sigma)))|^2}
\exp\left(-\int_\sigma^\tau\frac{\triangle g(\vec{x}(\phi(\rho)))}{|\nabla g(\vec{x}(\phi(\rho)))|^2}\, d\rho \right)\, d\sigma\,.
\end{aligned}
\end{equation}

\subsection{Time trace data}\label{subsec:titr}
The range condition becomes 
\begin{equation}\label{eqn:range_titr}
[\ul{u},\ol{u}]=h(0,T)\supseteq u_{act}(\Omega\times(0,T)) 
\end{equation}
and projection of the {\sc pde}  on the observation manifold $\{x_0\}\times(0,T)$ yields the iteration
\begin{equation}\label{eqn:proj_titr}
a_{k+1}'(h(t)) |\nabla u(x_0,t;a_k)|^2 + a_{k+1}(h(t)) \Delta u(x_0,t;a_k) = h'(t)-r(x_0,t) \,, \quad t\in(0,T)
\end{equation}
Here we again employ the {\sc ode} solution formula as above.
By abbreviating $u_k(x,t) = u(x,t;a_k)$ and 
using $a_{\rm act}(\tau_0)=a_{\rm act}(h(0))=\frac{b(x_0,0)-\gamma h(0)}{\partial_\nu u_0} =: a^0$
(or, if $\partial_\nu u_0=0$, assuming
$a_{\rm act}(\tau_0)=a_{\rm act}(h(0)) =: a^0$ to be known)
we obtain
\begin{equation}\label{eqn:a+_titr_PDE}
\begin{aligned}
a_{k+1}(\tau)=&\exp\left(-\int_{\tau_0}^\tau\frac{\triangle u_k(x_0,h^{-1}(\sigma))}{|\nabla u_k(x_0,h^{-1}(\sigma))|^2}\, d\sigma\right)a^0\\
&+ \int_{\tau_0}^\tau \frac{h'(h^{-1}(\sigma))-r(x_0,h^{-1}(\sigma))}{|\nabla u_k(x_0,h^{-1}(\tau))|^2}
\exp\left(-\int_\sigma^\tau \frac{\triangle u_k(x_0,h^{-1}(\rho))}{|\nabla u_k(x_0,h^{-1}(\rho))|^2}\, d\rho \right)\, d\sigma\,.
\end{aligned}
\end{equation}
Alternatively, with $x_0\in\partial\Omega$, we might slightly modify the iteration by inserting the nonlinear boundary conditions \eqref{eqn:bndy}, thus replacing the second term in \eqref{eqn:proj_titr} to obtain
\begin{equation}\label{eqn:proj_titr_alt}
a_{k+1}'(h(t)) |\nabla u_k(x_0,t)|^2 =
\frac{\gamma h(t)-b(x_0,t)}{\partial_\nu u_k(x_0,t)} \Delta u_k(x_0,t) + h'(t)-r(x_0,t) \,, \quad t\in(0,T).
\end{equation}
Then by inverting $h$ and integrating with respect to $\tau$ yields
\begin{equation}\label{eqn:a+_titr_PDE_bndy}
a_{k+1}(\tau)=a^0+\int_{\tau_0}^\tau\frac{\frac{\gamma \sigma-b(x_0,h^{-1}(\sigma))}{\partial_\nu u_k(x_0,h^{-1}(\sigma))} \Delta u_k(x_0,h^{-1}(\sigma)) + h'(h^{-1}(\sigma))-r(x_0,h^{-1}(\sigma))}{|\nabla u_k(x_0,h^{-1}(\sigma))|^2}\,d\sigma\,.
\end{equation}
A third scheme is obtained by relying exclusively on the nonlinear boundary conditions \eqref{eqn:bndy} and defining
\begin{equation}\label{eqn:a+_titr_bndy}
a_{k+1}(\tau)=\frac{b(x_0,h^{-1}(\tau))-\gamma \tau}{\partial_\nu u^D_k(x_0,h^{-1}(\tau))}
\end{equation}
where $u^D_k$ solves \eqref{eqn:u-nl} with the nonlinear impedance boundary conditions replaced by Dirichlet boundary conditions 
\begin{equation}\label{eqn:bndyD}
 u^D_k(x,t) = h(x,t), \quad x\in\Gamma,\quad t\in(0,T)
\end{equation}
on part of the boundary $\Gamma\subset\partial\Omega$, where we assume to have observations on $\Gamma\times(0,T)$. In case of $d>1$ space dimensions, $\Gamma$ should be a set of positive $d-1$ dimensional measure, whereas in 1-d it suffices to use $\Gamma=\{x_0\}$.

In all three cases, it is clearly crucial to have strict monotonicity of $h$
and boundedness away from zero of $|\nabla u_{act}|$, $\partial_\nu u_{act}$
(thus, by a perturbation argument, also of $|\nabla u_k|$, $\partial_\nu u_k$)
-- analogously to the assumption on $g$ in \eqref{eqn:Sigma}.

\section{Well-posedness of the forward problem \eqref{eqn:u-nl}, \eqref{eqn:bndy}}\label{sect:forwards}
For 
\begin{equation}\label{eq:PDEdata}
a\in C^{1+\alpha}(\R)\,, \quad a(s)\geq\underline{a}>0\,, \quad s\in\R\,, \quad r\in C^{\alpha,\alpha/2}(Q_T)\,, \quad b\in C^{\alpha,\alpha/2}(\partial\Omega\times(0,T))
\end{equation}
where $Q_T=\Omega\times(0,T)$ is the space-time cylinder, \cite[Theorem 6.1, p 453]{LadyzhenskajaSolonnikovUraltseva68} implies existence of a solution $u(a)\in H^{2+\alpha,1+\alpha/2}(Q_T)$ of \eqref{eqn:u-nl}, \eqref{eqn:bndy}, whose norm only depends on the norms of $a$, $r$, and $b$ in the mentioned spaces. To obtain Schauder estimates on $u=u(a)$, we first of all estimate the Schauder norm of the space- and time dependent coefficient $a(u)$
\[
\begin{aligned}
\|a(u)\|_{C^{\epsilon,\epsilon/2}(Q_T)}&
=\|a(u)\|_{C(Q_T)}+\sup_{(x,t)\not=(\tilde{x},\tilde{t})\in Q_T} \frac{a(u(x,t))-a(u(\tilde{x},\tilde{t}))}{\sqrt{|x-\tilde{x}|^2+|t-\tilde{t}|}^\epsilon}\\
&=\|a(u)\|_{C(Q_T)}+\sup_{(x,t)\not=(\tilde{x},\tilde{t})\in Q_T} \frac{a(u(x,t))-a(u(\tilde{x},\tilde{t}))}{|u(x,t)-u(\tilde{x},\tilde{t})|}
\frac{|u(x,t)-u(\tilde{x},\tilde{t})|}{\sqrt{|x-\tilde{x}|^2+|t-\tilde{t}|}^\epsilon}\\
&\leq\|a\|_{C(\R)}+ |a|_{C^{0+1}(\R)} |u|_{C^{\epsilon,\epsilon/2}(Q_T)}\\
&\leq\|a\|_{C(\R)}+ |a|_{C^{0+1}(\R)} C \|u\|_{H^{2+\alpha,1+\alpha/2}(Q_T)}\\
\end{aligned}
\]
for $0<\epsilon<2+\alpha-\frac{d+2}{2}$, by the embedding result \cite[Lemma 3.3, p 80]{LadyzhenskajaSolonnikovUraltseva68}.
Thus we can apply \cite[Theorem 6, page 65]{Friedman:1964} to the transformed version of \eqref{eqn:u-nl}, \eqref{eqn:bndy}
\[
\begin{aligned}
&D_t w-a(u)\triangle w= a(u)r  \mbox{ in }\Omega\times(0,T)\\ 
&w(0)=a(u_0) \mbox{ in }\Omega\\
&\partial_\nu w = b-\gamma u, \mbox{ on }\partial\Omega\times(0,T).
\end{aligned}
\]
Note the advantage of only needing $a(u)\in C^\epsilon$ in this non-divergence
form which allows us to apply the results from \cite{Friedman:1964}), where 
\[
w(x,t)= A(u(x,t))\quad A(s)=\int_0^s a(r)\, dr\,,
\]
to obtain
\[
\begin{aligned}
&\|w\|_{C^{\epsilon,\epsilon/2}(Q_T)}+\|D_t w\|_{C^{\epsilon,\epsilon/2}(Q_T)}
+\sum_{i=1}^n\|D_{x_i}w\|_{C^{\epsilon,\epsilon/2}(Q_T)}
+\sum_{i,j=1}^n\|D_{x_ix_j}w\|_{C^{\epsilon,\epsilon/2}(Q_T)}\\
&\qquad\qquad\leq K \Bigl(\|a(u_0)\|_{C^{2+\epsilon}(\Omega)} + \|a(u)r\|_{C^{\epsilon,\epsilon/2}(Q_T)}
+\|b-\gamma u\|_{C^{\epsilon,\epsilon/2}(\partial\Omega\times(0,T))}
\Bigr)\,.
\end{aligned}
\]
%Note that the first term on the right hand side yields a too high norm of $a$ unless $u_0$ is constant, so for simplicity we will assume $u_0=0$ from now on. 
The 
%BK 2020-10-16 second term
second and third terms 
on the right hand side can be estimated by $\|a(u)\|_{C^{\epsilon,\epsilon/2}(Q_T)}\, \|r\|_{C^{\epsilon,\epsilon/2}(Q_T)}+\|b\|_{C^{\epsilon,\epsilon/2}(\partial\Omega\times(0,T))}+\gamma C \|u\|_{H^{2+\alpha,1+\alpha/2}(Q_T)}$, which for $\epsilon\leq\alpha$ is already covered by the above assumptions and estimates.
From this we get 
\begin{equation}\label{eqn:estuC1}
\begin{aligned}
\|u\|_{C^{1,1/2}(Q_T)}
&\leq \|u\|_{C(Q_T)}+\sup_{(x,t)\not=(\tilde{x},\tilde{t})\in Q_T} \frac{u(x,t)-u(\tilde{x},\tilde{t})}{\sqrt{|x-\tilde{x}|^2+|t-\tilde{t}|}}\\
&\leq \|u\|_{C(Q_T)}+\sum_{i=1}^n\|D_{x_i}u\|_{C(Q_T)}+\|D_tu\|_{C(Q_T)}\\
&= \|u\|_{C(Q_T)}+\sum_{i=1}^n\|\frac{1}{a(u)}D_{x_i}w\|_{C(Q_T)}+\|\frac{1}{a(u)}D_tw\|_{C(Q_T)}\\
&\leq (1+\tfrac{1}{\underline{a}}) K \Bigl(
%BK 2020-10-16 add u0 term!
\|a(u_0)\|_{C^{2+\epsilon}(\Omega)}+
\|a(u)r\|_{C^{\epsilon,\epsilon/2}(Q_T)}\\
&\qquad+\|b\|_{C^{\epsilon,\epsilon/2}(\partial\Omega\times(0,T))}+\gamma C \|u\|_{H^{2+\alpha,1+\alpha/2}(Q_T)}\Bigr)
\ \leq 
%BK 2020-10-16 C_{\rm arb}\,.
C_{{\rm arb}\, u_0}\,.
\end{aligned}
\end{equation}
where 
%BK 2020-10-16 $C_{\rm arb}$
$C_{{\rm arb}\, u_0}$
denotes a generic constant depending only on the quantities 
\begin{equation}\label{eq:Carb}
\underline{a},\ \|a\|_{C^\alpha(\R)},\ \|r\|_{C^{\alpha,\alpha/2}(Q_T)},\ \|b\|_{C^{\alpha,\alpha/2}(\partial\Omega\times(0,T)}, \
%BK 2020-10-16 add u0 term!
\|a(u_0)\|_{C^{2+\alpha}(\Omega)}.
\end{equation}
Thus the above estimate on $a(u)$ can be improved to
\begin{equation}\label{eqn:estauCalpha}
\begin{aligned}
\|a(u)\|_{C^{\alpha,\alpha/2}(Q_T)}
&=\|a(u)\|_{C(Q_T)}+\sup_{(x,t)\not=(\tilde{x},\tilde{t})\in Q_T} \frac{a(u(x,t))-a(u(\tilde{x},\tilde{t}))}{|u(x,t)-u(\tilde{x},\tilde{t})|^\alpha}
\Bigl(\frac{|u(x,t)-u(\tilde{x},\tilde{t})|}{\sqrt{|x-\tilde{x}|^2+|t-\tilde{t}|}}\Bigr)^\alpha\\
&\leq\|a\|_{C^\alpha(\R)} \|u\|_{C^{1,1/2}(Q_T)}^\alpha
\ \leq 
%BK 2020-10-16 C_{\rm arb}\,.
C_{{\rm arb}\, u_0}\,.
\end{aligned}
\end{equation}
and we can apply the above result from \cite{Friedman:1964} with $\alpha$ in place of $\epsilon$ to obtain 
\[
\begin{aligned}
&\|w\|_{C^{\alpha,\alpha/2}(Q_T)}+\|D_t w\|_{C^{\alpha,\alpha/2}(Q_T)}
+\sum_{i=1}^n\|D_{x_i}w\|_{C^{\alpha,\alpha/2}(Q_T)}
+\sum_{i,j=1}^n\|D_{x_ix_j}w\|_{C^{\alpha,\alpha/2}(Q_T)}\\
&\quad \leq K \Bigl(\|a(u_0)\|_{C^{2+\alpha}(\Omega)} +\|a(u)\|_{C^{\alpha,\alpha/2}(Q_T)}\, \|r\|_{C^{\alpha,\alpha/2}(Q_T)}
+\|b\|_{C^{\alpha,\alpha/2}(\partial\Omega\times(0,T))}+\gamma 
\underbrace{\|u\|_{C^{\alpha,\alpha/2}(Q_T)}}_{\leq \|u\|_{C^{1,1/2}(Q_T)}}
\Bigr)
\end{aligned}
\]
In particular, we get, besides \eqref{eqn:estuC1}
\begin{equation}\label{eqn:estuxCalpha}
\|\nabla u\|_{C^{\alpha,\alpha/2}(Q_T)} = \|\frac{1}{a(u)} \nabla w\|_{C^{\alpha,\alpha/2}(Q_T)}
\ \leq C_{{\rm arb}\, u_0}
\end{equation}
and
\begin{equation}\label{eqn:estuxxCalpha}
\|\triangle u\|_{C^{\alpha,\alpha/2}(Q_T)} = \|\frac{1}{a(u)} (\triangle w-\frac{a'(u)}{a(u)^2}|\nabla w|^2)\|_{C^{\alpha,\alpha/2}(Q_T)}
\ \leq C_{{\rm arb}\, u_0}
\end{equation}

\begin{proposition}
For $a$, $r$, $b$ as in \eqref{eq:PDEdata} and $u_0$ such that $a(u_0)\in C^{2+\alpha}(\Omega)$, there exists a unique solution $u$ of \eqref{eqn:u-nl}, \eqref{eqn:bndy}, which satisfies 
\[
\|a(u)D_tu\|_{C^{\alpha,\alpha/2}(Q_T)} + \|\nabla u\|_{C^{\alpha,\alpha/2}(Q_T)} + \|\triangle u\|_{C^{\alpha,\alpha/2}(Q_T)} \leq C_{{\rm arb}\, u_0}
\]
with some constant $C_{{\rm arb}\, u_0}$ depending only on the quantities in \eqref{eq:Carb}. 
%BK 2020-10-16 and $\|a(u_0)\|_{C^{2+\alpha}(\Omega)}$.
\end{proposition}

\begin{remark}
Note that in case of constant initial data $u_0$, 
%BK 2020-10-16
condition
$a(u_0)\in C^{2+\alpha}(\Omega)$ does not imply any additional assumptions on the smoothness of $a$.
\end{remark}

\section{Contractivity with final time data}\label{sec:fiti_contractivity}

We start with some estimates on the iteration errors in $a$ and $u$. Our goal is to prove that 
\[
\|a_{k+1}-a_{\rm act}\|\leq q\|a_k-a_{\rm act}\|
\] for some $q\in(0,1)$ and some appropriate norm $\|\cdot\|$, that is, the fixed point iteration defined in Section \oldref{subsec:fiti} is contractive.

\smallskip

Abbreviating $\hat{a}=a_k-a_{\rm act}$, $\hat{a}_+=a_{k+1}-a_{\rm act}$, $\hat{u}(x,t)=u(x,t;a_k)-u(x,t;a_{\rm act})=u_k(x,t)-u_{\rm act}(x,t)$ we get the following error representations.

In one spatial dimension \eqref{eqn:a+1-d} yields
\[
\begin{aligned}
\hat{a}_+(g(x))&=\tfrac{1}{g_x(x)}\int_0^x D_t \hat{u}(\xi,T)\, d\xi \\
\hat{a}_+'(g(x))g_x(x)&= \tfrac{1}{g_x(x)} D_t \hat{u}(x,T) + (\tfrac{1}{g_x})_x(x)\int_0^x D_t \hat{u}(\xi,T)\, d\xi
\end{aligned}
\]
In higher space dimensions, we get from \eqref{eqn:a+h-d} and \eqref{eqn:ODEa+}
\[
\begin{aligned}
\hat{a}_+(\tau)&=\int_0^\tau \frac{D_t \hat{u}(\vec{x}(\phi(\sigma)),T;a_k)}{|\nabla g(\vec{x}(\phi(\sigma)))|^2}
\exp\left(-\int_\sigma^\tau\frac{\triangle g(\vec{x}(\phi(\rho)))}{|\nabla g(\vec{x}(\phi(\rho)))|^2}\, d\rho \right)\, d\sigma
\\
\hat{a}_+'(\tau)&=\frac{1}{|\nabla g(\vec{x}(\phi(\tau)))|^2} \left( -\hat{a}_+(\tau)\triangle g(\vec{x}(\phi(\tau))) + D_t \hat{u}(\vec{x}(\phi(\tau)),T)\right)\,.
\end{aligned}
\]
Hence, there exists a constant $C_{g\Sigma}=C(\kappa,\|g\|_{C^{2+\alpha}(\Omega)},\|\vec{x}\|_{C^{0+1}(0,1)})$ depending only on the data $g$ and the observation manifold cf. \eqref{eqn:Sigma}, \eqref{eqn:phi} such that   
\begin{equation}\label{eqn:ahatDtuhat}
\|a_{k+1}-a_{\rm act}\|_{C^{2+\alpha}(J)}= \|\hat{a}_+\|_{C^{2+\alpha}(J)} 
\leq C_{g\Sigma} \|D_t\hat{u}(T)\|_{C^{1+\alpha}(\Omega)}\,.
\end{equation}

Thus it remains to estimate the norm of $z(T)=D_t\hat{u}(T)$ by a small multiple of $\hat{a}$. To this end, observe that %$\tilde{u}_{\rm act}=D_t u_{\rm act}$ 
$\hat{u}=u-u_{\rm act}$ and $z=D_t\hat{u}$ solve
\begin{equation}\label{eqn:PDEuhat}
\begin{aligned}
&D_t \hat{u} -\triangle (\bar{a} \hat{u})
= \triangle \Bigl(\int_0^{u}\hat{a}(s)\, ds\Bigr) = \nabla\cdot(\hat{a}(u)\nabla u)  \mbox{ in }\Omega\times(0,T)\\ 
&\hat{u}(0)=0 \mbox{ in }\Omega\\
&\partial_\nu (\bar{a} \hat{u})+\gamma \hat{u}=-\hat{a}(u)\partial_\nu u, \mbox{ on }\partial\Omega\times(0,T)
\end{aligned}
\end{equation}
and
\begin{equation}\label{eqn:PDEz}
\begin{aligned}
&D_t z -\triangle (\bar{a} z)
= \triangle \Bigl(\hat{a}(u)(D_t u_{\rm act}+z) +(\bar{a}_{1,0} D_t u_{\rm act} + \bar{a}_{1,1} z) \hat{u}\Bigr)
  \mbox{ in }\Omega\times(0,T)\\ 
&z(0)=D_t\hat{u}(0)= \nabla\cdot(\hat{a}(u_0)\nabla u_0) \mbox{ in }\Omega\\
&\partial_\nu (\bar{a} z)+\gamma z = -\partial_\nu (\hat{a}(u)(D_t u_{\rm act}+z) + (\bar{a}_{1,0} D_t u_{\rm act} + \bar{a}_{1,1} z) \hat{u}\Bigr), \mbox{ on }\partial\Omega\times(0,T),
\end{aligned}
\end{equation}
where
\begin{equation}\label{eqn:aij}
\bar{a} = \bar{a}_{0,0} = \int_0^1 a_{\rm act}(u_{\rm act}+\theta\hat{u})\, d\theta\,, \qquad
\bar{a}_{i,j} = \int_0^1 a_{\rm act}^{(i)}(u_{\rm act}+\theta\hat{u})\theta^j\, d\theta
\end{equation}
so that $D_t \bar{a} = \bar{a}_{1,0} D_t u_{\rm act} + \bar{a}_{1,1} z$.
Let us assume now that $a_{\rm act}$ is positive and bounded away from zero by some constant $\underline{a}$
\[
a_{\rm act}(s)\geq\underline{a}>0
\]
and we have Dirichlet boundary conditions in equation~ \eqref{eqn:bndy},
so that the boundary conditions in \eqref{eqn:PDEz} simplify to homogeneous
Dirichlet ones.
Similarly to 
%BK 2020-10-16 \cite{KaltenbacherRundell:2019d},
\cite{KaltenbacherRundell:2020a},
we can make use of the exponential decay of $z$ and $D_t u_{\rm act}$ in order
to achieve contractivity.
To this end we multiply \eqref{eqn:PDEz} with $e^{\mu t}$ which for $z_\mu(x,t)=e^{\mu t} z(x,t)$, $\tilde{u}_\mu(x,t)=e^{\mu t} D_tu_{\rm act}(x,t)$  
yields 
\begin{equation}\label{eqn:PDEzmu}
\begin{aligned}
&D_t z_\mu - \triangle(\bar{a} z_\mu) = \triangle[y_1 z_\mu + y_2 \tilde{u}_\mu]  \mbox{ in }\Omega\times(0,T)\\ 
&z_\mu(0)=\nabla\cdot(\hat{a}(u_0)\nabla u_0) \mbox{ in }\Omega\\
&z_\mu = 0 \mbox{ on }\partial\Omega\times(0,T)
\end{aligned}
\end{equation}
with the multipliers
\begin{equation}\label{eqn:y1y2}
y_1 = \bar{a}_{1,1} \hat{u} +  \hat{a}(u)\,, \qquad
y_2 = \bar{a}_{1,0} \hat{u} + \hat{a}(u) \,,
\end{equation}
where we choose $\mu\in(0,\underline{a}\lambda_1)$ with $\lambda_1>0$ the smallest eigenvalue of $-\triangle$ (with Dirichlet boundary conditions) so that $-\triangle(\bar{a} \cdot)-\mu I$ is still elliptic.
Then maximal parabolic regularity yields 
\begin{equation}\label{eqn:CAz}
\begin{aligned}
&\|z_\mu\|_{W^{1,p}(0,t;L^p(\Omega))}+\|z_\mu\|_{L^p(0,t;W^{2,p}(\Omega))} \\
&\qquad\leq C_{\underline{a},\mu,p} \Bigl(\|\triangle[y_1 z_\mu + y_2 \tilde{u}_\mu]\|_{L^p(\Omega\times(0,t))} 
+\|\nabla\cdot(\hat{a}(u_0)\nabla u_0)\|_{W^{2(1-1/p),p}(\Omega)}\Bigr)
\end{aligned}
\end{equation}
with a constant $C_{\underline{a},\mu,p}$ independent of $t$,
(cf., e.g., \cite[Proposition 3.1]{Latushkin2006}, \cite[Theorem 2.1]{DenkHieberPruess2007}).
Here, on one hand, we can use continuity of the embeddings $W^{\theta,p}(0,t)\to C(0,t)$ and $W^{2-2\theta,p}(\Omega)\to C^{1,\alpha}(\Omega)$ for $\theta\in(0,1)$, $\theta>\frac{1}{p}$, $1-2\theta>\frac{d}{p}+\alpha$ (which can always be achieved by choosing $p\in [1,\infty)$ sufficiently large) and apply interpolation and the 
the fact that the norm of the embedding $W^{\theta,p}(0,t)\to C(0,t)$ is independent of $t$ (by Morrey's inequality) to obtain
\begin{equation}\label{eqn:estz}
\begin{aligned}
\|e^{\mu t} z(t)\|_{C^{1+\alpha}(\Omega)} & \leq \|z_\mu\|_{C(0,t;C^{1+\alpha}(\Omega))}
\ \leq C_{W^{\theta,p},C}^{\mathbb{R^+}} C_{W^{2-2\theta,p},C^{1+\alpha}}^\Omega \|z_\mu\|_{W^{\theta,p}(0,t;W^{2-2\theta,p}(\Omega))}\\
&\leq C_{W^{\theta,p},C}^{\mathbb{R^+}} C_{W^{2-2\theta,p},C^{1+\alpha}}^\Omega(\|z_\mu\|_{W^{1,p}(0,t;L^p(\Omega))}+\|z_\mu\|_{L^p(0,t;W^{2,p}(\Omega))}).
\end{aligned}
\end{equation}
On the other hand, we can estimate the right hand side in \eqref{eqn:CAz} by
\[
\|\triangle[y_1 z_\mu + y_2 \tilde{u}_\mu]\|_{L^p(\Omega\times(0,t))} \leq
\sum_{i=1}^2\|\triangle[y_i w_i]\|_{L^p(\Omega\times(0,t))} 
\]
where $w_1=z_\mu$, $w_2=\tilde{u}_\mu$ and both terms can be estimated in the same way:
\begin{equation}\label{eqn:yiwi}
\begin{aligned}
&\|\triangle[y_i w_i]\|_{L^p(\Omega\times(0,t))} \ = \|\triangle y_i \, w_i + 2 \nabla y_i \cdot \nabla w_i + y_i \triangle w_i\|_{L^p(\Omega\times(0,t))}\\
&\leq \|\triangle y_i\|_{L^p(0,t;L^p(\Omega))} \|w_i\|_{L^\infty(0,t;L^\infty(\Omega))} 
+ 2 \|\nabla y_i\|_{L^\infty(0,t;L^{qp/(q-p)}(\Omega)} \|\nabla w_i\|_{L^p(0,t;L^q(\Omega)} \\
&\qquad+ \|y_i\|_{L^\infty(\Omega\times(0,t))} \|\triangle w_i\|_{L^p(0,t;L^p(\Omega))}\\
&\leq \Bigl(C_{W^{\theta,p},C}^{\mathbb{R^+}} C_{W^{2-2\theta,p},C}^\Omega \|\triangle y_i\|_{L^p(0,t;L^p(\Omega))}
+ 2 C_{W^{2,p},W^{1,q}}^\Omega \|\nabla y_i\|_{L^\infty(0,t;L^{qp/(q-p)}(\Omega)}
+ \|y_i\|_{L^\infty(\Omega\times(0,t))}
\Bigr)\\
&\qquad\qquad \cdot(\|w_i\|_{W^{1,p}(0,t;L^p(\Omega))}+\|w_i\|_{L^p(0,t;W^{2,p}(\Omega))})
\end{aligned}
\end{equation}
with 
\begin{equation}\label{eqn:pq}
2-\frac{d}{p}>0\,, \quad 1-\frac{d}{p}>-\frac{d}{q}\,, \quad q>p\,.
\end{equation}

The $y_i$- terms in \eqref{eqn:yiwi} can be estimated by (see \eqref{eqn:y1y2}, \eqref{eqn:aij})
%y_1 = \bar{a}_{1,1} \hat{u} + \hat{a}(u)\,, \qquad
%y_2 = \bar{a}_{1,0} \hat{u} + \hat{a}(u) \,. 
\[
\begin{aligned}
\|y_1\|_{L^\infty(\Omega\times(0,t))}
&\leq c_a \|\hat{u}\|_{L^\infty(\Omega\times(0,t))} + \|\hat{a}\|_{C(J)}
\\[1ex]
\|\nabla y_1\|_{L^\infty(0,t;L^{qp/(q-p)}(\Omega)}
&=\|(\bar{a}_{2,1}\hat{u}+\hat{a}'(u))\nabla u_{\rm act} + (\bar{a}_{1,1}+\bar{a}_{2,2}\hat{u}+\hat{a}'(u))\nabla \hat{u} 
\|_{L^\infty(0,t;L^{qp/(q-p)}(\Omega)}\\
&\leq |\Omega|^{(q-p)/(qp)} \Bigl(
(\tfrac12 c_a \|\hat{u}\|_{C(\Omega\times(0,t))} + \|\hat{a}'\|_{C(J)}) \|\nabla u_{\rm act}\|_{C(\Omega\times(0,t))}\\
&\qquad\qquad
+ ((\tfrac12+\tfrac13\|\hat{u}\|_{C(\Omega\times(0,t))})c_a  + \|\hat{a}'\|_{C(J)}) \|\nabla \hat{u}\|_{L^\infty(0,t;L^{qp/(q-p)}(\Omega)}\Bigr)
\end{aligned}
\]
\[
\begin{aligned}
&\|\triangle y_1\|_{L^p(0,t;L^p(\Omega))}
=\|(\bar{a}_{2,1}\hat{u}+\hat{a}'(u))\triangle u_{\rm act} + (\bar{a}_{1,1}+\bar{a}_{2,2}\hat{u}+\hat{a}'(u))\triangle \hat{u}\\
&\qquad\qquad\qquad\qquad +
(\bar{a}_{3,1}\hat{u}+\hat{a}''(u))|\nabla u_{\rm act}|^2 + 
2(\bar{a}_{2,1}+\bar{a}_{3,2}\hat{u}+\hat{a}''(u))\nabla u_{\rm act}\cdot\nabla \hat{u}\\
&\qquad\qquad\qquad\qquad +
(2\bar{a}_{2,2}+\bar{a}_{3,3}\hat{u}+\hat{a}''(u))|\nabla \hat{u}|^2\|_{L^p(0,t;L^p(\Omega))}\\
&\qquad\leq 
\Bigl(\tfrac12 c_a \|\hat{u}\|_{C(\Omega\times(0,t))} + \|\hat{a}'\|_{C(J)} + 
C_{W^{2,p},W^{1,2p}}^\Omega \bigl((\tfrac12 + \tfrac56\|\hat{u}\|_{C(\Omega\times(0,t))}) c_a + 2\|\hat{a}''\|_{C(J)})\bigr)
\Bigr) \\
&\hspace*{5cm}\cdot\|u_{\rm act}\|_{L^p(0,t;W^{2,p}(\Omega))}\\
&\qquad\quad
+\Bigl((\tfrac12+\tfrac13\|\hat{u}\|_{C(\Omega\times(0,t))})c_a  + \|\hat{a}'\|_{C(J)} + 
C_{W^{2,p},W^{1,2p}}^\Omega \bigl((\tfrac76 + \tfrac{7}{12}\|\hat{u}\|_{C(\Omega\times(0,t))}) c_a + 2\|\hat{a}''\|_{C(J)})\bigr)
\Bigr) \\
&\hspace*{5cm}\cdot\|\hat{u}\|_{L^p(0,t;W^{2,p}(\Omega))},
\end{aligned}
\]
where we have used continuity of the embedding $W^{2,p}(\Omega)\to W^{1,2p}(\Omega)$ (cf. \eqref{eqn:pq}) and
the fact that with $c_a:=\|a_{\rm act}\|_{C^3(J)}$ we have 
$\|\bar{a}_{i,j}\|_{L^\infty(\Omega\times(0,t))}\leq \frac{1}{j+1} c_a$.
We can estimate $y_2$ analogously.

%For $\hat{u}$ we get from \eqref{eqn:PDEuhat} that 
%\[
%\|\hat{u}\|_{C^{2+\epsilon,1+\epsilon/2}(Q_T)}\leq 
%K\|\nabla\cdot(\hat{a}(u)\nabla u)\|_{C^{\epsilon,\epsilon/2}(Q_T)}\leq C_{arb} \|\hat{a}\|_{C^{1+\epsilon}(J)}
%\]

%For $\hat{u}$ we get from \eqref{eqn:PDEuhat} (but with homogeneous Dirichlet boundary conditions) and the same maximal parabolic regularity estimate as above that 
%\[
%\begin{aligned}
%&\|\hat{u}\|_{L^{2p}(0,t;W^{1,2p}(\Omega))}+\|\hat{u}\|_{L^p(0,t;W^{2,p}(\Omega))}
%\leq C_{\underline{a},\mu,p} \|\nabla\cdot(\hat{a}(u)\nabla u)\|_{L^p(\Omega\times(0,t))} \\
%&\leq C_{\underline{a},\mu,p} \Bigl(\|\hat{a}\|_{C(J)}
%\Bigl(\|u_{\rm act}\|_{L^p(0,t;W^{2,p}(\Omega))}+\|\hat{u}\|_{L^p(0,t;W^{2,p}(\Omega))}\Bigr)\\
%&\qquad\qquad+ \|\hat{a}'\|_{C(J)} \Bigl(\|u_{\rm act}\|_{L^{2p}(0,t;W^{1,2p}(\Omega))}+\|\hat{u}\|_{L^{2p}(0,t;W^{1,2p}(\Omega))}\Bigr)^2\Bigr)\,,
%\end{aligned}
%\]
%hence for $\|\hat{a}\|_{C^1(J)}\leq\frac{1}{2C_{\underline{a},\mu,p}}$
%\[
%\|\hat{u}\|_{L^p(0,t;W^{2,p}(\Omega))}
%\leq 2 C_{\underline{a},\mu,p} \|\hat{a}\|_{C^1(J)} \|u_{\rm act}\|_{L^p(0,t;W^{2,p}(\Omega))}\,.
%\]

For $\hat{u}$  defined in \eqref{eqn:PDEuhat} (but with homogeneous Dirichlet
boundary conditions) 
%BK 2020-10-16 and
with
the same maximal parabolic regularity as well as
embedding and interpolation estimate as above, we get that 
\[
\begin{aligned}
\|\hat{u}\|_{C(\Omega\times(0,t))}+\|\hat{u}\|_{L^p(0,t;W^{2,p}(\Omega))}
&\leq C_{\underline{a},p} \|\nabla\cdot(\hat{a}(u)\nabla u)\|_{L^p(\Omega\times(0,t))} \\
&\leq C_{\underline{a},p} 
\Bigl(\|\hat{a}\|_{C(J)}\|u\|_{L^p(0,t;W^{2,p}(\Omega))}
+ \|\hat{a}'\|_{C(J)} \|u\|_{L^{2p}(0,t;W^{1,2p}(\Omega))}^2\Bigr)\\
&\leq C_{\underline{a},p} \Bigl(C_{rdu_0} + \bigl(C_{W^{1,p}(L^p)\cap L^p(W^{2,p}),L^{2p}(W^{1,2p})}^{Q_T} C_{rdu_0}\bigr)^2\Bigr)
\|\hat{a}\|_{C^1(J)} \,,
\end{aligned}
\]
where 
\[
C_{rdu_0} = C_{\underline{a},p} (\|r\|_{L^p(\Omega\times(0,t))}+\|d\|_{W_{p,tr}} +\|u_0\|_{W^{2(1-1/p),p}(\Omega)})
\]
with $W_{p,tr} = W^{1-1/(2p),p}(0,t;L_p(\partial\Omega))\cap L^p(0,t;W^{2-1/p}(\Omega))$.

For $w_2=\tilde{u}_\mu = e^{\mu t}D_t u_{\rm act}$ we get, with the same maximal parabolic regularity result and the fact that it solves (note that we are in the Dirichlet boundary setting $\gamma=\infty$ now)
\begin{equation}\label{eqn:PDEutilmu}
\begin{aligned}
&D_t \tilde{u}_\mu  -\triangle[a_{\rm act}(u_{\rm act})\tilde{u}_\mu] -\mu \tilde{u}_\mu= \tilde{r}_\mu   \mbox{ in }\Omega\times(0,T)\\ 
&\tilde{u}_\mu(0)=\nabla\cdot(a_{\rm act}(u_0)\nabla u_0)+D_tr(0) \mbox{ in }\Omega\\
&\tilde{u}_\mu = D_td \mbox{ on }\partial\Omega\times(0,T)
\end{aligned}
\end{equation}
the estimate
\begin{equation}\label{eqn:CAutil}
\begin{aligned}
&\|\tilde{u}_\mu\|_{W^{1,p}(0,t;L^p(\Omega))}+\|\tilde{u}_\mu\|_{L^p(0,t;W^{2,p}(\Omega))} \\
&\qquad\leq C_{\underline{a},\mu,p} \Bigl(\|\tilde{r}_\mu\|_{L^p(\Omega\times(0,t))} + \|\tilde{d}_\mu\|_{W_{p,tr}}
%\\&\qquad\qquad
+  \|\nabla\cdot(a_{\rm act}(u_0)\nabla u_0)+D_tr(0)\|_{W^{2(1-1/p),p}(\Omega)}\Bigr)\\
\end{aligned}
\end{equation}
where $\tilde{r}_\mu(x,t)=e^{\mu t} D_tr(x,t)$, $\tilde{d}_\mu(x,t)=e^{\mu t} D_td(x,t)$.

The initial data terms in \eqref{eqn:CAz}, \eqref{eqn:CAutil}
can be estimated by 
\[
\|\nabla\cdot(\tilde{a}(u_0)\nabla u_0)\|_{W^{2(1-1/p),p}(\Omega)}\leq C \|\tilde{a}\|_{C^{2+\alpha}(J)} \|u_0\|_{W^{4-2/p,p}(\Omega)}
\]
for $\tilde{a}\in\{\hat{a},a_{\rm act}\}$, $\alpha\geq 1-2/p$.

Note that for spatial dimension $d=3$,
according to \eqref{eqn:pq} we have to choose $p>\frac32$
and therefore are in the the regime $1-\frac{1}{2p}>\frac{1}{p}$ where (cf. 
%BK 2020-10-16 DenkHieberPruess2007) 
\cite{DenkHieberPruess2007}) 
compatibility conditions 
\begin{equation}\label{eqn:compat}
\nabla\cdot(\hat{a}(u_0)\nabla u_0) = 0\,, \quad 
\nabla\cdot(a_{\rm act}(u_0)\nabla u_0)+D_tr(0) = D_t d(0)
\mbox{ on }\partial\Omega
\end{equation}
need to be assumed.

From the above we get an estimate of the form
\[
\begin{aligned}
&\|z_\mu\|_{W^{1,p}(0,t;L^p(\Omega))}+\|z_\mu\|_{L^p(0,t;W^{2,p}(\Omega))}\\ 
&\leq C \|\hat{a}\|_{C^{2+\alpha}(J)} \Bigl(\|z_\mu\|_{W^{1,p}(0,t;L^p(\Omega))}+\|z_\mu\|_{L^p(0,t;W^{2,p}(\Omega))}
\\&\qquad \qquad
+\|u_0\|_{W^{4-2/p,p}(\Omega)} + \|D_tr(0)\|_{W^{2(1-1/p),p}(\Omega)}
+ \|\tilde{r}_\mu\|_{L^p(\Omega\times(0,t))} + \|\tilde{d}_\mu\|_{W_{p,tr}}\Bigr)
\end{aligned}
\]
with some constant $C$ independent of $t$, hence, for $\|\hat{a}\|_{C^{2+\alpha}(J)}\leq \frac{1}{2C}$ 
\begin{equation}\label{eqn:estz2}
\begin{aligned}
&\|z_\mu\|_{W^{1,p}(0,t;L^p(\Omega))}+\|z_\mu\|_{L^p(0,t;W^{2,p}(\Omega))}\\ 
&\leq C \|\hat{a}\|_{C^{2+\alpha}(J)} \Bigl(\|u_0\|_{W^{4-2/p,p}(\Omega)} + \|D_tr(0)\|_{W^{2(1-1/p),p}(\Omega)}
+ \|\tilde{r}_\mu\|_{L^p(\Omega\times(0,t))} + \|\tilde{d}_\mu\|_{W_{p,tr}}\Bigr)\,.
\end{aligned}
\end{equation}

Thus, provided that for some $p>\frac{d}{2}$, (cf. \eqref{eqn:pq})
$D_t r$, $D_t d$  decay exponentially
\begin{equation}\label{eqn:Dtrdecay}
\begin{aligned}
&\|D_tr(t)\|_{L^p(\Omega)}\leq C_r e^{-\mu_r t} \,, \quad  
\|D_td(t)\|_{W^{2-1/p}(\Omega)}\leq C_d e^{-\mu_d t} \,, \quad t>0 \,, \\
&\|e^{\mu_d\cdot}D_td\|_{W^{1-1/(2p),p}(0,\infty;L_p(\partial\Omega))}:=D_2 <\infty\,,
\end{aligned}
\end{equation} 
for some $C_r, C_d, \mu_r, \mu_d >0$ so that
\begin{equation}
\begin{aligned}
&\|\tilde{r}_\mu\|_{L^p(\Omega\times(0,\infty))}
\leq (\frac{C_r}{p(\mu_r-\mu)})^{1/p}=:R<\infty\\
&\|\tilde{d}_\mu\|_{L^p(0,t;W^{2-1/p}(\Omega))}\leq
(\frac{C_d}{p(\mu_d-\mu)})^{1/p}=:D_1<\infty\\
\end{aligned}
\end{equation} 
for any $\mu\in(0,\min\{\mu_r,\mu_d,\underline{a}\lambda_1\})\!\!$
\footnote{Note that the characterisation of general Sobolev spaces via fractional derivatives seems to be still open and so characterisation of the last condition in \eqref{eqn:Dtrdecay} by exponential decay of some fractional derivative of $D_t d$ is not possible} we get from  
\eqref{eqn:CAz}, \eqref{eqn:estz2}
\[
\|z(T)\|_{C^{1+\alpha}(\Omega)}  \leq C e^{-\mu tT} 
\Bigl(\|u_0\|_{W^{4-2/p,p}(\Omega)} + \|D_tr(0)\|_{W^{2(1-1/p),p}(\Omega)}
+ R + D_1 + D_2\Bigr) \ \|\hat{a}\|_{C^{2+\alpha}(J)}
\]
with $C$ independent of $T$.

In view of \eqref{eqn:ahatDtuhat} we have thus proven the following contractivity result.

\begin{theorem}\label{th:contr_fiti}
Assume $g\in C^{3+\alpha}(\Omega)$, \eqref{eqn:Sigma} with $\vec{x}\in C^2([0,1])$, 
$a_{\rm act}\in C^3(\R)$, $a_{\rm act}=a_{le}(\tau)$ for $\tau<\ul{u}$, $a_{\rm act}(\tau)=a_{ri}(\tau)$ for $\tau>\ol{u}$, 
$a_{\rm act}(\tau)\geq \underline{a}>0$, $\tau\in\mathbb{R}$,
$u_0\in W^{4-2/p,p}(\Omega)\cap C^{2+\alpha}(\Omega)$,
$D_tr(0)\in W^{2(1-1/p),p}(\Omega)$
\eqref{eqn:Dtrdecay}, for some $p>\frac{d}{2}$, $\alpha\geq 1-2/p$, 
%BK 2020-10-16
and, in case $d=3$, \eqref{eqn:compat}.
Then for $T$ large enough and $\|a_0-a_{\rm act}\|_{C^{2+\alpha}(J)}$ small enough there exists a constant $q\in(0,1)$ such that for all $k\in\N$
\[
\|a_{k+1}-a_{\rm act}\|_{C^{2+\alpha}(J)} \leq q \|a_k-a_{\rm act}\|_{C^{2+\alpha}(J)}\,.
\]
\end{theorem}

\begin{corollary}
Under the assumptions of Theorem \oldref{th:contr_fiti}, $a$ is uniquely determined.
\end{corollary}

\begin{remark}
In reality, instead of $g$ we have noisy data lying only in $L^p(\Omega)$, with a noise level $\delta$ that is also only given (if at all) with respect to the $L^p$ norm. Smoothing this data leads to an approximate version of $g$ in $C^{2+\alpha}$ with a (typically larger) noise level $\tilde{\delta}$ in the $C^{2+\alpha}$ norm. Due to contractivity of the scheme, this noise will propagate boundedly through the iteration so that after sufficiently many -- $O(\log(1/\tilde{\delta})$ -- iterations (but without having to stop early, in principle)
we end up with an $O(\tilde{\delta})$ accuracy in the reconstruction.
For details, see, e.g., \cite[Section 4.5]{KaltenbacherRundell:2019b}.
\end{remark}

\section{Reconstructions}\label{sect:recons}

In this section we will show the results of numerical experiments to recover 
the function $a(u)$  with the two different types of data measurements.
The computations will be done in one space dimension.

The first of these is when we are only able to obtain ``census-type''
information and thus measure $g(x) := u(x,T)$ for some
fixed time $T$ and $x\in \Sigma$ where
$\Sigma$ is a curve in $\Omega$ whose endpoints lie
at different points on $\partial\Omega$.
Since our reconstructions will be based on one space variable
we simply use $u(x,T)$ as the data and ensure that our imposed initial/boundary
conditions allows the range condition to hold.

The second is when we are able to measure 
%BK 2020-10-16 multiple types of data from a
%BK 2020-10-16 single experimental run: specifically both the final data and a
%BK 2020-10-16 boundary measurement of either the solution $u$ or its normal derivative at a
the solution $u$ at a 
fixed point $x_0\in\partial\Omega$.
The range condition 
%BK 2020-10-16 the
then 
can be imposed by driving that boundary $x=x_0$
where we measure $u(x_0,t)$ with sufficiently large supplied heat flux.

To procure data for the reconstructions a direct solver based on a
Crank-Nicolson scheme produced output values and data
values were produced from these by sampling at a relatively small
number $N_x$ and $N_t$ of points in both the spatial and temporal directions:
This sampled data was then both interpolated to a full working size to obtain
data values $g(x)$ and $h(t)$ commensurate to the grid being used by the solver
used in the inverse problem and such that the value of $g(x)$ was
smoothed by an $H^2$ filter while that of $h(t)$ by an $H^1$ filter.

Figure~\oldref{fig:titi_titi} shows reconstructions of a function $a(u)$
that goes well beyond a low-degree polynomial from time trace data,
%BK 2020-10-16
using the iteration \eqref{eqn:a+_titr_bndy}, that clearly outperformed 
the two schemes \eqref{eqn:a+_titr_PDE}, \eqref{eqn:a+_titr_PDE_bndy} 
in all our experiments.
The initial approximation was $a(u) = 1$.
The leftmost figure shows the reconstructions of selected iterations;
the final one shown was the effective numerical convergence.
The noise level here was $0.01\%$.
The rightmost figure shows the convergence of the $n^{\rm th}$ iterate
%BK 2020-10-16 $a_n(u)$; 
$a_n$;
specifically $\|a_n - a_{act}\|$ for both $L^2$ and $L^\infty$.
As suggested by the left figure, the convergence is initially quite
rapid but slows down considerably.
Certainly, at the scale of the graph the tenth and final iteration
would be indistinguishable.

%  \PiCTeX file for BB4
\input colordvi

\input pictex
\font\smallsymbol = cmmi8
\newdimen\xfiglen \newdimen\yfiglen
\xfiglen=2.6 true in
\yfiglen=1.8 true in
\newbox\figurelegendone
\newbox\figurelegendtwo
\newbox\figurelegendthree
\newbox\figurelegendfour
\newbox\figurelegendfive
\newbox\figureone
\newbox\figuretwo
\newbox\figurethree 
\newbox\figurefour
\newbox\figurefive
\newbox\figuresix
%%%%%%%%%%%%%%%%%%%

\setbox\figurelegendone=\hbox{
\beginpicture
  \setcoordinatesystem units <0.3\xfiglen,0.37\yfiglen> %point at 0 -0.7
  \setplotarea x from 0 to 0.8, y from 0 to 0.8
\linethickness=0.7pt
\footnotesize
\ninerm
  \Red{\relax
  \putrule from 0 0.0 to 0.25 0.0 }\relax
  \put {iteration 25} [l] at 0.35 0.0
  \Orange{\relax
  \setdashes <2pt>
  \putrule from 0 0.2 to 0.25 0.2 }\relax
  \put {iteration 5} [l] at 0.35 0.2
  \OliveGreen{\relax
  \setdashes <2pt>
   \putrule from 0 0.4 to 0.25 0.4 }\relax
  \put {iteration 2} [l] at 0.35 0.4
  \Blue{\relax
   \setsolid
   \putrule from 0 0.6 to 0.25 0.6 }\relax
  \put {iteration 1} [l] at 0.35 0.6
\setsolid
\setplotsymbol ({\sixrm .})
  \putrule from 0 0.8 to 0.25 0.8
  \put {actual $a(u)$} [l] at 0.35 0.8
\endpicture
\relax
}
\setbox\figurelegendtwo=\hbox{
\footnotesize
\beginpicture
  \setcoordinatesystem units <0.2\xfiglen,0.3\yfiglen> %point at 0 -0.7
  \setplotarea x from 0 to 0.8, y from 0 to 0.6
  \put {${\tiny{\circ}}$ $\|a_n-a_{act}\|_{L^\infty}$} [l] at 0 0.25
  \put {${\tiny{\bullet}}$ $\|a_n- a_{act}\|_{L^2}$} [l] at 0 0.6
\endpicture
\relax
}
\setbox\figurelegendthree=\hbox{
\footnotesize
\beginpicture
  \setcoordinatesystem units <0.3\xfiglen,0.4\yfiglen> %point at 0 -0.7
  \setplotarea x from 0 to 0.8, y from 0 to 0.65
\linethickness=0.7pt
\ninerm
  \Red{\relax
  \putrule from 0 0.0 to 0.25 0.0 }\relax
  \put {iteration 5} [l] at 0.35 0.0
  \putrule from 0 0.2 to 0.25 0.2  %\relax
  \OliveGreen{\relax
  \putrule from 0 0.2 to 0.25 0.2 }\relax
  \put {iteration 2} [l] at 0.35 0.2
  \Blue{\relax
   \putrule from 0 0.4 to 0.25 0.4 }\relax
  \putrule from 0 0.4 to 0.2 0.4
  \put {iteration 1} [l] at 0.35 0.4
\setsolid
\setplotsymbol ({\sixrm .})
  \putrule from 0 0.6 to 0.25 0.6
  \put {actual $a(u)$} [l] at 0.35 0.6
\endpicture
\relax
}
\setbox\figurelegendfour=\hbox{
\footnotesize
\beginpicture
  \setcoordinatesystem units <0.3\xfiglen,0.4\yfiglen> %point at 0 -0.7
  \setplotarea x from 0 to 0.8, y from 0 to 0.65
\linethickness=0.7pt
\ninerm
  %\Red{\relax
  %\putrule from 0 0.0 to 0.25 0.0 }\relax
  %\put {iteration 5} [l] at 0.35 0.0
  %\putrule from 0 0.2 to 0.25 0.2  %\relax
  \Red{\relax
  \putrule from 0 0.2 to 0.25 0.2 }\relax
  \put {5\% noise} [l] at 0.35 0.2
  \Green{\relax
   \putrule from 0 0.4 to 0.25 0.4 }\relax
  \put {1\% noise} [l] at 0.35 0.4
\setsolid
\setplotsymbol ({\sixrm .})
  \putrule from 0 0.6 to 0.25 0.6
  \put {actual $a(u)$} [l] at 0.35 0.6
\endpicture
\relax
}
\setbox\figurelegendfive=\hbox{
\footnotesize
\beginpicture
  \setcoordinatesystem units <0.3\xfiglen,0.4\yfiglen> %point at 0 -0.7
  \setplotarea x from 0 to 0.8, y from 0.2 to 0.65
\linethickness=0.7pt
\ninerm
  \Red{\relax
  \putrule from 0 0.2 to 0.25 0.2  }\relax
  \put {iteration 3} [l] at 0.35 0.2
  \Blue{\relax
   \putrule from 0 0.4 to 0.25 0.4 }\relax
  \put {iteration 1} [l] at 0.35 0.4
\setsolid
\setplotsymbol ({\sixrm .})
  \putrule from 0 0.6 to 0.25 0.6
  \put {actual $a(u)$} [l] at 0.35 0.6
\endpicture
\relax
}
\setbox\figureone=\vbox{\hsize=\xfiglen
\beginpicture
\eightrm
  \setcoordinatesystem units <0.45\xfiglen,0.52\yfiglen>  point at 0.0 0.5
  \setplotarea x from 0 to 2.2, y from 0.5 to 2.5
  \axis bottom shiftedto y=0.5 ticks short numbered from 0 to 2.2 by 0.5 /
  \axis left ticks short numbered from 0.5 to 2.5 by 0.5 /
\footnotesize
\put {$u$} [rb] at 2.17 0.54
\put {$a(u)$} [lt] at 0.025 2.5
\put {\copy\figurelegendone} [rt] at 2.2 2.53
\setsolid
\setplotsymbol ({\sixrm .})
\Black{\relax  % actual $a(u) = (0.75 + exp(-50*(u-1).^2) + 2*u.*exp(-u) + 0.5*sin(3.5*u))$
\plot
         0    0.7500
    0.0210    0.8279
    0.0420    0.9040
    0.0631    0.9779
    0.0841    1.0497
    0.1051    1.1191
    0.1261    1.1860
    0.1472    1.2504
    0.1682    1.3119
    0.1892    1.3706
    0.2102    1.4264
    0.2313    1.4790
    0.2523    1.5284
    0.2733    1.5745
    0.2943    1.6173
    0.3154    1.6566
    0.3364    1.6924
    0.3574    1.7246
    0.3784    1.7533
    0.3995    1.7784
    0.4205    1.7998
    0.4415    1.8177
    0.4625    1.8319
    0.4835    1.8426
    0.5046    1.8498
    0.5256    1.8535
    0.5466    1.8539
    0.5676    1.8510
    0.5887    1.8450
    0.6097    1.8360
    0.6307    1.8245
    0.6517    1.8107
    0.6728    1.7954
    0.6938    1.7797
    0.7148    1.7651
    0.7358    1.7538
    0.7569    1.7486
    0.7779    1.7530
    0.7989    1.7704
    0.8199    1.8042
    0.8410    1.8562
    0.8620    1.9260
    0.8830    2.0102
    0.9040    2.1018
    0.9250    2.1907
    0.9461    2.2649
    0.9671    2.3122
    0.9881    2.3229
    1.0091    2.2913
    1.0302    2.2172
    1.0512    2.1060
    1.0722    1.9675
    1.0932    1.8140
    1.1143    1.6580
    1.1353    1.5105
    1.1563    1.3790
    1.1773    1.2680
    1.1984    1.1785
    1.2194    1.1092
    1.2404    1.0572
    1.2614    1.0193
    1.2825    0.9924
    1.3035    0.9736
    1.3245    0.9611
    1.3455    0.9533
    1.3665    0.9494
    1.3876    0.9486
    1.4086    0.9508
    1.4296    0.9557
    1.4506    0.9631
    1.4717    0.9729
    1.4927    0.9851
    1.5137    0.9996
    1.5347    1.0163
    1.5558    1.0350
    1.5768    1.0556
    1.5978    1.0780
    1.6188    1.1021
    1.6399    1.1276
    1.6609    1.1546
    1.6819    1.1826
    1.7029    1.2117
    1.7240    1.2416
    1.7450    1.2721
    1.7660    1.3030
    1.7870    1.3342
    1.8080    1.3654
    1.8291    1.3965
    1.8501    1.4272
    1.8711    1.4574
    1.8921    1.4869
    1.9132    1.5154
    1.9342    1.5429
    1.9552    1.5691
    1.9762    1.5938
    1.9973    1.6169
    2.0183    1.6383
    2.0393    1.6578
    2.0603    1.6752
    2.0814    1.6905
    2.1024    1.7035
 /\relax}\relax

\setplotsymbol ({\fiverm .})
\setsolid
\Blue{\relax  % iteration 1
\plot
         0    1.0330
    0.0210    0.9433
    0.0420    0.9225
    0.0631    0.9449
    0.0841    0.9790
    0.1051    1.0107
    0.1261    1.0386
    0.1472    1.0642
    0.1682    1.0886
    0.1892    1.1116
    0.2102    1.1327
    0.2313    1.1523
    0.2523    1.1708
    0.2733    1.1884
    0.2943    1.2048
    0.3154    1.2200
    0.3364    1.2339
    0.3574    1.2468
    0.3784    1.2588
    0.3995    1.2697
    0.4205    1.2797
    0.4415    1.2888
    0.4625    1.2971
    0.4835    1.3045
    0.5046    1.3109
    0.5256    1.3166
    0.5466    1.3215
    0.5676    1.3257
    0.5887    1.3293
    0.6097    1.3321
    0.6307    1.3342
    0.6517    1.3357
    0.6728    1.3367
    0.6938    1.3371
    0.7148    1.3372
    0.7358    1.3371
    0.7569    1.3371
    0.7779    1.3377
    0.7989    1.3395
    0.8199    1.3433
    0.8410    1.3498
    0.8620    1.3595
    0.8830    1.3724
    0.9040    1.3887
    0.9250    1.4080
    0.9461    1.4286
    0.9671    1.4487
    0.9881    1.4669
    1.0091    1.4817
    1.0302    1.4915
    1.0512    1.4956
    1.0722    1.4937
    1.0932    1.4858
    1.1143    1.4722
    1.1353    1.4540
    1.1563    1.4317
    1.1773    1.4069
    1.1984    1.3815
    1.2194    1.3571
    1.2404    1.3335
    1.2614    1.3118
    1.2825    1.2920
    1.3035    1.2743
    1.3245    1.2601
    1.3455    1.2486
    1.3665    1.2382
    1.3876    1.2292
    1.4086    1.2222
    1.4296    1.2166
    1.4506    1.2120
    1.4717    1.2085
    1.4927    1.2063
    1.5137    1.2050
    1.5347    1.2048
    1.5558    1.2055
    1.5768    1.2072
    1.5978    1.2098
    1.6188    1.2131
    1.6399    1.2171
    1.6609    1.2217
    1.6819    1.2269
    1.7029    1.2329
    1.7240    1.2394
    1.7450    1.2462
    1.7660    1.2532
    1.7870    1.2606
    1.8080    1.2684
    1.8291    1.2765
    1.8501    1.2846
    1.8711    1.2928
    1.8921    1.3010
    1.9132    1.3093
    1.9342    1.3174
    1.9552    1.3256
    1.9762    1.3337
    1.9973    1.3414
    2.0183    1.3492
    2.0393    1.3567
    2.0603    1.3635
    2.0814    1.3701
    2.1024    1.3768
 /\relax}\relax

\setdashes <2pt>
\OliveGreen{\relax    % iteration 3
\plot
         0    1.0173
    0.0210    0.9051
    0.0420    0.8736
    0.0631    0.9073
    0.0841    0.9688
    0.1051    1.0331
    0.1261    1.0924
    0.1472    1.1472
    0.1682    1.1990
    0.1892    1.2477
    0.2102    1.2921
    0.2313    1.3337
    0.2523    1.3733
    0.2733    1.4115
    0.2943    1.4477
    0.3154    1.4810
    0.3364    1.5115
    0.3574    1.5398
    0.3784    1.5660
    0.3995    1.5899
    0.4205    1.6113
    0.4415    1.6309
    0.4625    1.6484
    0.4835    1.6633
    0.5046    1.6758
    0.5256    1.6863
    0.5466    1.6948
    0.5676    1.7015
    0.5887    1.7064
    0.6097    1.7094
    0.6307    1.7106
    0.6517    1.7101
    0.6728    1.7082
    0.6938    1.7049
    0.7148    1.7005
    0.7358    1.6954
    0.7569    1.6907
    0.7779    1.6879
    0.7989    1.6892
    0.8199    1.6968
    0.8410    1.7127
    0.8620    1.7382
    0.8830    1.7739
    0.9040    1.8204
    0.9250    1.8773
    0.9461    1.9386
    0.9671    1.9983
    0.9881    2.0502
    1.0091    2.0891
    1.0302    2.1097
    1.0512    2.1101
    1.0722    2.0901
    1.0932    2.0503
    1.1143    1.9935
    1.1353    1.9242
    1.1563    1.8449
    1.1773    1.7618
    1.1984    1.6812
    1.2194    1.6084
    1.2404    1.5409
    1.2614    1.4807
    1.2825    1.4271
    1.3035    1.3803
    1.3245    1.3447
    1.3455    1.3172
    1.3665    1.2936
    1.3876    1.2731
    1.4086    1.2577
    1.4296    1.2462
    1.4506    1.2371
    1.4717    1.2307
    1.4927    1.2270
    1.5137    1.2257
    1.5347    1.2268
    1.5558    1.2303
    1.5768    1.2360
    1.5978    1.2437
    1.6188    1.2530
    1.6399    1.2635
    1.6609    1.2754
    1.6819    1.2890
    1.7029    1.3043
    1.7240    1.3207
    1.7450    1.3376
    1.7660    1.3547
    1.7870    1.3725
    1.8080    1.3914
    1.8291    1.4111
    1.8501    1.4310
    1.8711    1.4505
    1.8921    1.4698
    1.9132    1.4893
    1.9342    1.5085
    1.9552    1.5276
    1.9762    1.5462
    1.9973    1.5639
    2.0183    1.5816
    2.0393    1.5983
    2.0603    1.6130
    2.0814    1.6272
    2.1024    1.6414
 /\relax}\relax

\Orange{\relax  % iteration 5
\plot
         0    1.0155
    0.0210    0.8974
    0.0420    0.8629
    0.0631    0.8994
    0.0841    0.9695
    0.1051    1.0457
    0.1261    1.1179
    0.1472    1.1851
    0.1682    1.2482
    0.1892    1.3067
    0.2102    1.3595
    0.2313    1.4086
    0.2523    1.4554
    0.2733    1.5008
    0.2943    1.5441
    0.3154    1.5838
    0.3364    1.6201
    0.3574    1.6535
    0.3784    1.6842
    0.3995    1.7120
    0.4205    1.7367
    0.4415    1.7590
    0.4625    1.7785
    0.4835    1.7946
    0.5046    1.8075
    0.5256    1.8177
    0.5466    1.8251
    0.5676    1.8301
    0.5887    1.8328
    0.6097    1.8330
    0.6307    1.8310
    0.6517    1.8268
    0.6728    1.8206
    0.6938    1.8125
    0.7148    1.8026
    0.7358    1.7918
    0.7569    1.7815
    0.7779    1.7739
    0.7989    1.7722
    0.8199    1.7798
    0.8410    1.7994
    0.8620    1.8328
    0.8830    1.8809
    0.9040    1.9450
    0.9250    2.0242
    0.9461    2.1097
    0.9671    2.1915
    0.9881    2.2600
    1.0091    2.3067
    1.0302    2.3244
    1.0512    2.3114
    1.0722    2.2689
    1.0932    2.1988
    1.1143    2.1060
    1.1353    1.9986
    1.1563    1.8805
    1.1773    1.7611
    1.1984    1.6491
    1.2194    1.5515
    1.2404    1.4638
    1.2614    1.3872
    1.2825    1.3200
    1.3035    1.2622
    1.3245    1.2193
    1.3455    1.1876
    1.3665    1.1615
    1.3876    1.1395
    1.4086    1.1235
    1.4296    1.1122
    1.4506    1.1039
    1.4717    1.0988
    1.4927    1.0969
    1.5137    1.0977
    1.5347    1.1012
    1.5558    1.1077
    1.5768    1.1170
    1.5978    1.1286
    1.6188    1.1418
    1.6399    1.1564
    1.6609    1.1724
    1.6819    1.1905
    1.7029    1.2108
    1.7240    1.2323
    1.7450    1.2542
    1.7660    1.2762
    1.7870    1.2989
    1.8080    1.3229
    1.8291    1.3480
    1.8501    1.3732
    1.8711    1.3978
    1.8921    1.4221
    1.9132    1.4464
    1.9342    1.4703
    1.9552    1.4940
    1.9762    1.5170
    1.9973    1.5388
    2.0183    1.5605
    2.0393    1.5806
    2.0603    1.5981
    2.0814    1.6148
    2.1024    1.6315
 /\relax}\relax

\setsolid
\Red{\relax   % iteration 25
\plot
         0    1.0150
    0.0210    0.8950
    0.0420    0.8598
    0.0631    0.8980
    0.0841    0.9725
    0.1051    1.0552
    0.1261    1.1352
    0.1472    1.2100
    0.1682    1.2796
    0.1892    1.3429
    0.2102    1.3989
    0.2313    1.4501
    0.2523    1.4988
    0.2733    1.5462
    0.2943    1.5918
    0.3154    1.6338
    0.3364    1.6720
    0.3574    1.7070
    0.3784    1.7388
    0.3995    1.7671
    0.4205    1.7918
    0.4415    1.8135
    0.4625    1.8317
    0.4835    1.8460
    0.5046    1.8564
    0.5256    1.8634
    0.5466    1.8670
    0.5676    1.8676
    0.5887    1.8655
    0.6097    1.8610
    0.6307    1.8544
    0.6517    1.8455
    0.6728    1.8341
    0.6938    1.8199
    0.7148    1.8028
    0.7358    1.7834
    0.7569    1.7639
    0.7779    1.7478
    0.7989    1.7404
    0.8199    1.7470
    0.8410    1.7722
    0.8620    1.8187
    0.8830    1.8879
    0.9040    1.9808
    0.9250    2.0944
    0.9461    2.2131
    0.9671    2.3187
    0.9881    2.3945
    1.0091    2.4293
    1.0302    2.4172
    1.0512    2.3612
    1.0722    2.2685
    1.0932    2.1444
    1.1143    1.9967
    1.1353    1.8380
    1.1563    1.6736
    1.1773    1.5155
    1.1984    1.3744
    1.2194    1.2580
    1.2404    1.1595
    1.2614    1.0775
    1.2825    1.0085
    1.3035    0.9510
    1.3245    0.9107
    1.3455    0.8843
    1.3665    0.8667
    1.3876    0.8556
    1.4086    0.8514
    1.4296    0.8531
    1.4506    0.8598
    1.4717    0.8700
    1.4927    0.8831
    1.5137    0.8987
    1.5347    0.9162
    1.5558    0.9377
    1.5768    0.9629
    1.5978    0.9903
    1.6188    1.0184
    1.6399    1.0462
    1.6609    1.0746
    1.6819    1.1051
    1.7029    1.1382
    1.7240    1.1722
    1.7450    1.2056
    1.7660    1.2381
    1.7870    1.2707
    1.8080    1.3049
    1.8291    1.3403
    1.8501    1.3755
    1.8711    1.4094
    1.8921    1.4422
    1.9132    1.4746
    1.9342    1.5059
    1.9552    1.5367
    1.9762    1.5662
    1.9973    1.5938
    2.0183    1.6204
    2.0393    1.6441
    2.0603    1.6639
    2.0814    1.6818
    2.1024    1.6996
 /\relax}\relax

\endpicture
}
\yfiglen=1.5 true in
\setbox\figuretwo=\vbox{\hsize=\xfiglen
% This is for the $a$ reconstruction
\beginpicture
\eightrm
  \setcoordinatesystem units <0.04\xfiglen,2.1\yfiglen>  point at 0.0 0
  \setplotarea x from 0 to 25, y from 0.0 to 0.6
  \axis bottom shiftedto y=0 ticks short numbered from 0 to 25 by 5 /
  \axis left ticks short numbered from 0.0 to 0.6 by 0.1 /
\put {\copy\figurelegendtwo} [rt] at 24 0.58
\footnotesize
\put {$n$} [b] at 25 0.01
\multiput {${\tiny{\circ}}$} at   % L^2 norm differences
    0    0.4083
    1    0.2513
    2    0.1733
    3    0.1398
    4    0.1213
    5    0.1086
    6    0.0995
    7    0.0932
    8    0.0887
    9    0.0853
   10    0.0825
   11    0.0801
   12    0.0781
   13    0.0764
   14    0.0749
   15    0.0737
   16    0.0726
   17    0.0717
   18    0.0710
   19    0.0703
   20    0.0698
   21    0.0693
   22    0.0690
   23    0.0687
   24    0.0684
   25    0.0682
/
\multiput {${\tiny{\bullet}}$} at  % L_inf norm differences
    0    0.5695
    1    0.3722
    2    0.2279
    3    0.2171
    4    0.2221
    5    0.2163
    6    0.2064
    7    0.1957
    8    0.1863
    9    0.1782
   10    0.1716
   11    0.1666
   12    0.1625
   13    0.1592
   14    0.1566
   15    0.1545
   16    0.1530
   17    0.1517
   18    0.1507
   19    0.1498
   20    0.1491
   21    0.1485
   22    0.1480
   23    0.1476
   24    0.1472
   25    0.1470
/
\endpicture
}
%%%%%
%
\yfiglen=1.8 true in
\setbox\figurethree=\vbox{\hsize=\xfiglen
\beginpicture
\eightrm
  \setcoordinatesystem units <0.45\xfiglen,0.5\yfiglen>  point at 0.0 0.5
  \setplotarea x from 0 to 2.2, y from 0.5 to 2.5
  \axis bottom shiftedto y=0.5 ticks short numbered from 0 to 2.2 by 0.5 /
  \axis left ticks short numbered from 0.5 to 2.5 by 0.5 /
\footnotesize
\put {$u$} [rb] at 2.17 0.54
\put {$a(u)$} [lt] at 0.025 2.5
\put {\copy\figurelegendthree} [rt] at 2.15 2.5
\setsolid
\setplotsymbol ({\sixrm .})

\setsolid
\setplotsymbol ({\fiverm .})
\Blue{\relax  % iteration 1
\plot
    0.0000    0.5695
    0.0210    0.6193
    0.0420    0.6727
    0.0631    0.7274
    0.0841    0.7835
    0.1051    0.8419
    0.1261    0.9004
    0.1471    0.9568
    0.1681    1.0097
    0.1892    1.0587
    0.2102    1.1042
    0.2312    1.1469
    0.2522    1.1877
    0.2732    1.2270
    0.2943    1.2639
    0.3153    1.2984
    0.3363    1.3300
    0.3573    1.3589
    0.3783    1.3852
    0.3994    1.4092
    0.4204    1.4306
    0.4414    1.4494
    0.4624    1.4659
    0.4834    1.4799
    0.5044    1.4912
    0.5255    1.4998
    0.5465    1.5058
    0.5675    1.5094
    0.5885    1.5104
    0.6095    1.5089
    0.6306    1.5052
    0.6516    1.4995
    0.6726    1.4924
    0.6936    1.4849
    0.7146    1.4782
    0.7357    1.4742
    0.7567    1.4754
    0.7777    1.4849
    0.7987    1.5056
    0.8197    1.5408
    0.8407    1.5919
    0.8618    1.6595
    0.8828    1.7404
    0.9038    1.8279
    0.9248    1.9133
    0.9458    1.9875
    0.9669    2.0415
    0.9879    2.0678
    1.0089    2.0610
    1.0299    2.0188
    1.0509    1.9426
    1.0720    1.8384
    1.0930    1.7161
    1.1140    1.5875
    1.1350    1.4606
    1.1560    1.3419
    1.1770    1.2370
    1.1981    1.1491
    1.2191    1.0767
    1.2401    1.0200
    1.2611    0.9767
    1.2821    0.9457
    1.3032    0.9242
    1.3242    0.9098
    1.3452    0.9013
    1.3662    0.8971
    1.3872    0.8962
    1.4083    0.8981
    1.4293    0.9025
    1.4503    0.9092
    1.4713    0.9178
    1.4923    0.9284
    1.5133    0.9411
    1.5344    0.9564
    1.5554    0.9742
    1.5764    0.9935
    1.5974    1.0146
    1.6184    1.0381
    1.6395    1.0642
    1.6605    1.0918
    1.6815    1.1205
    1.7025    1.1493
    1.7235    1.1787
    1.7445    1.2094
    1.7656    1.2407
    1.7866    1.2728
    1.8076    1.3050
    1.8286    1.3377
    1.8496    1.3704
    1.8707    1.4036
    1.8917    1.4365
    1.9127    1.4683
    1.9337    1.4994
    1.9547    1.5293
    1.9758    1.5577
    1.9968    1.5853
    2.0178    1.6115
    2.0388    1.6359
    2.0598    1.6585
    2.0808    1.6803
    2.1019    1.7018
 /\relax}\relax

\OliveGreen{\relax    % iteration 2
\plot
    0.0000    0.6792
    0.0210    0.7387
    0.0420    0.8024
    0.0631    0.8676
    0.0841    0.9346
    0.1051    1.0043
    0.1261    1.0740
    0.1471    1.1412
    0.1681    1.2042
    0.1892    1.2624
    0.2102    1.3162
    0.2312    1.3666
    0.2522    1.4144
    0.2732    1.4603
    0.2943    1.5030
    0.3153    1.5426
    0.3363    1.5784
    0.3573    1.6108
    0.3783    1.6397
    0.3994    1.6655
    0.4204    1.6880
    0.4414    1.7071
    0.4624    1.7231
    0.4834    1.7358
    0.5044    1.7452
    0.5255    1.7513
    0.5465    1.7542
    0.5675    1.7540
    0.5885    1.7508
    0.6095    1.7445
    0.6306    1.7355
    0.6516    1.7243
    0.6726    1.7116
    0.6936    1.6983
    0.7146    1.6858
    0.7357    1.6763
    0.7567    1.6723
    0.7777    1.6772
    0.7987    1.6942
    0.8197    1.7264
    0.8407    1.7754
    0.8618    1.8415
    0.8828    1.9207
    0.9038    2.0059
    0.9248    2.0879
    0.9458    2.1570
    0.9669    2.2041
    0.9879    2.2220
    1.0089    2.2055
    1.0299    2.1531
    1.0509    2.0667
    1.0720    1.9528
    1.0930    1.8216
    1.1140    1.6852
    1.1350    1.5513
    1.1560    1.4265
    1.1770    1.3162
    1.1981    1.2239
    1.2191    1.1478
    1.2401    1.0881
    1.2611    1.0425
    1.2821    1.0098
    1.3032    0.9871
    1.3242    0.9719
    1.3452    0.9628
    1.3662    0.9582
    1.3872    0.9572
    1.4083    0.9590
    1.4293    0.9633
    1.4503    0.9700
    1.4713    0.9787
    1.4923    0.9893
    1.5133    1.0022
    1.5344    1.0176
    1.5554    1.0355
    1.5764    1.0549
    1.5974    1.0761
    1.6184    1.0996
    1.6395    1.1256
    1.6605    1.1530
    1.6815    1.1814
    1.7025    1.2098
    1.7235    1.2385
    1.7445    1.2684
    1.7656    1.2986
    1.7866    1.3293
    1.8076    1.3599
    1.8286    1.3906
    1.8496    1.4212
    1.8707    1.4519
    1.8917    1.4819
    1.9127    1.5105
    1.9337    1.5381
    1.9547    1.5643
    1.9758    1.5886
    1.9968    1.6117
    2.0178    1.6332
    2.0388    1.6526
    2.0598    1.6700
    2.0808    1.6865
    2.1019    1.7025
 /\relax}\relax

\Red{\relax  % iteration 25
\plot
    0.0000    0.7239
    0.0210    0.7873
    0.0420    0.8552
    0.0631    0.9246
    0.0841    0.9960
    0.1051    1.0701
    0.1261    1.1442
    0.1471    1.2156
    0.1681    1.2824
    0.1892    1.3441
    0.2102    1.4010
    0.2312    1.4542
    0.2522    1.5047
    0.2732    1.5528
    0.2943    1.5976
    0.3153    1.6389
    0.3363    1.6761
    0.3573    1.7096
    0.3783    1.7392
    0.3994    1.7655
    0.4204    1.7881
    0.4414    1.8070
    0.4624    1.8224
    0.4834    1.8344
    0.5044    1.8427
    0.5255    1.8475
    0.5465    1.8488
    0.5675    1.8467
    0.5885    1.8413
    0.6095    1.8326
    0.6306    1.8210
    0.6516    1.8071
    0.6726    1.7917
    0.6936    1.7759
    0.7146    1.7612
    0.7357    1.7497
    0.7567    1.7445
    0.7777    1.7486
    0.7987    1.7657
    0.8197    1.7988
    0.8407    1.8495
    0.8618    1.9178
    0.8828    1.9996
    0.9038    2.0873
    0.9248    2.1709
    0.9458    2.2403
    0.9669    2.2858
    0.9879    2.2996
    1.0089    2.2766
    1.0299    2.2151
    1.0509    2.1178
    1.0720    1.9919
    1.0930    1.8490
    1.1140    1.7021
    1.1350    1.5596
    1.1560    1.4282
    1.1770    1.3131
    1.1981    1.2174
    1.2191    1.1392
    1.2401    1.0782
    1.2611    1.0318
    1.2821    0.9986
    1.3032    0.9757
    1.3242    0.9604
    1.3452    0.9512
    1.3662    0.9468
    1.3872    0.9458
    1.4083    0.9478
    1.4293    0.9524
    1.4503    0.9594
    1.4713    0.9684
    1.4923    0.9794
    1.5133    0.9927
    1.5344    1.0085
    1.5554    1.0268
    1.5764    1.0467
    1.5974    1.0684
    1.6184    1.0923
    1.6395    1.1189
    1.6605    1.1468
    1.6815    1.1757
    1.7025    1.2045
    1.7235    1.2338
    1.7445    1.2641
    1.7656    1.2947
    1.7866    1.3259
    1.8076    1.3569
    1.8286    1.3880
    1.8496    1.4189
    1.8707    1.4499
    1.8917    1.4802
    1.9127    1.5090
    1.9337    1.5369
    1.9547    1.5632
    1.9758    1.5877
    1.9968    1.6109
    2.0178    1.6325
    2.0388    1.6521
    2.0598    1.6697
    2.0808    1.6863
    2.1019    1.7025
 /\relax}\relax

\endpicture
}
\yfiglen=1.8 true in
\setbox\figurefour=\vbox{\hsize=\xfiglen
\beginpicture
\eightrm
  \setcoordinatesystem units <0.45\xfiglen,0.5\yfiglen>  point at 0.0 0.5
  \setplotarea x from 0 to 2.2, y from 0.5 to 2.5
  \axis bottom shiftedto y=0.5 ticks short numbered from 0 to 2.2 by 0.5 /
  \axis left ticks short numbered from 0.5 to 2.5 by 0.5 /
\put {\copy\figurelegendfour} [rt] at 2.15 2.5
\footnotesize
\put {$u$} [rb] at 2.17 0.54
\put {$a(u)$} [lt] at 0.025 2.5
\setsolid
\setplotsymbol ({\sixrm .})

\setsolid
\setplotsymbol ({\fiverm .})
\Green{\relax  % iteration 25
\plot
    0.0052    0.9500
    0.0262    0.9900
    0.0472    1.0221
    0.0682    1.0518
    0.0892    1.0841
    0.1102    1.1196
    0.1312    1.1589
    0.1521    1.2019
    0.1731    1.2486
    0.1941    1.2985
    0.2151    1.3512
    0.2361    1.4051
    0.2571    1.4590
    0.2781    1.5116
    0.2991    1.5617
    0.3201    1.6083
    0.3411    1.6505
    0.3621    1.6879
    0.3831    1.7202
    0.4040    1.7478
    0.4250    1.7721
    0.4460    1.7947
    0.4670    1.8165
    0.4880    1.8369
    0.5090    1.8541
    0.5300    1.8660
    0.5510    1.8709
    0.5720    1.8681
    0.5930    1.8580
    0.6140    1.8418
    0.6350    1.8212
    0.6559    1.7987
    0.6769    1.7773
    0.6979    1.7605
    0.7189    1.7517
    0.7399    1.7543
    0.7609    1.7711
    0.7819    1.8027
    0.8029    1.8481
    0.8239    1.9035
    0.8449    1.9643
    0.8659    2.0250
    0.8869    2.0804
    0.9079    2.1268
    0.9288    2.1624
    0.9498    2.1862
    0.9708    2.1962
    0.9918    2.1890
    1.0128    2.1595
    1.0338    2.1030
    1.0548    2.0169
    1.0758    1.9045
    1.0968    1.7747
    1.1178    1.6383
    1.1388    1.5063
    1.1598    1.3862
    1.1807    1.2817
    1.2017    1.1949
    1.2227    1.1259
    1.2437    1.0726
    1.2647    1.0325
    1.2857    1.0033
    1.3067    0.9823
    1.3277    0.9678
    1.3487    0.9580
    1.3697    0.9520
    1.3907    0.9489
    1.4117    0.9480
    1.4326    0.9491
    1.4536    0.9520
    1.4746    0.9568
    1.4956    0.9635
    1.5166    0.9724
    1.5376    0.9835
    1.5586    0.9969
    1.5796    1.0129
    1.6006    1.0318
    1.6216    1.0538
    1.6426    1.0792
    1.6636    1.1075
    1.6846    1.1387
    1.7055    1.1723
    1.7265    1.2075
    1.7475    1.2437
    1.7685    1.2805
    1.7895    1.3171
    1.8105    1.3530
    1.8315    1.3877
    1.8525    1.4209
    1.8735    1.4525
    1.8945    1.4824
    1.9155    1.5112
    1.9365    1.5396
    1.9574    1.5684
    1.9784    1.5986
    1.9994    1.6305
    2.0204    1.6644
    2.0414    1.7003
    2.0624    1.7383
    2.0834    1.7777
    2.1044    1.8177
 /\relax}\relax

\Red{\relax  % iteration 25
\plot
    0.0108    1.1178
    0.0319    1.1310
    0.0530    1.1446
    0.0741    1.1591
    0.0952    1.1751
    0.1162    1.1932
    0.1373    1.2136
    0.1584    1.2365
    0.1795    1.2622
    0.2006    1.2904
    0.2217    1.3212
    0.2428    1.3544
    0.2638    1.3899
    0.2849    1.4274
    0.3060    1.4665
    0.3271    1.5066
    0.3482    1.5470
    0.3693    1.5871
    0.3904    1.6269
    0.4114    1.6665
    0.4325    1.7069
    0.4536    1.7489
    0.4747    1.7920
    0.4958    1.8345
    0.5169    1.8727
    0.5380    1.9026
    0.5590    1.9212
    0.5801    1.9272
    0.6012    1.9215
    0.6223    1.9066
    0.6434    1.8863
    0.6645    1.8654
    0.6856    1.8492
    0.7066    1.8429
    0.7277    1.8506
    0.7488    1.8752
    0.7699    1.9166
    0.7910    1.9713
    0.8121    2.0320
    0.8332    2.0893
    0.8542    2.1347
    0.8753    2.1625
    0.8964    2.1714
    0.9175    2.1647
    0.9386    2.1469
    0.9597    2.1202
    0.9808    2.0842
    1.0018    2.0355
    1.0229    1.9700
    1.0440    1.8844
    1.0651    1.7807
    1.0862    1.6655
    1.1073    1.5468
    1.1284    1.4320
    1.1494    1.3271
    1.1705    1.2355
    1.1916    1.1584
    1.2127    1.0959
    1.2338    1.0466
    1.2549    1.0089
    1.2760    0.9806
    1.2970    0.9597
    1.3181    0.9444
    1.3392    0.9334
    1.3603    0.9254
    1.3814    0.9195
    1.4025    0.9155
    1.4236    0.9129
    1.4446    0.9119
    1.4657    0.9125
    1.4868    0.9146
    1.5079    0.9186
    1.5290    0.9246
    1.5501    0.9326
    1.5712    0.9430
    1.5922    0.9558
    1.6133    0.9715
    1.6344    0.9906
    1.6555    1.0134
    1.6766    1.0400
    1.6977    1.0701
    1.7188    1.1038
    1.7399    1.1407
    1.7609    1.1805
    1.7820    1.2222
    1.8031    1.2652
    1.8242    1.3086
    1.8453    1.3523
    1.8664    1.3963
    1.8875    1.4399
    1.9085    1.4829
    1.9296    1.5258
    1.9507    1.5690
    1.9718    1.6134
    1.9929    1.6601
    2.0140    1.7099
    2.0351    1.7629
    2.0561    1.8190
    2.0772    1.8777
    2.0983    1.9380
    2.1194    1.9990
 /\relax}\relax

\endpicture
}
\yfiglen=1.8 true in
\setbox\figurefive=\vbox{\hsize=\xfiglen
\beginpicture
\eightrm
  \setcoordinatesystem units <0.35\xfiglen,0.5\yfiglen>  point at 0.0 0.0
  \setplotarea x from 0 to 3, y from 0.0 to 2.2
  \axis bottom shiftedto y=0.0 ticks short numbered from 0 to 3 by 1 /
  \axis left ticks short numbered from 0.0 to 2.2 by 0.5 /
\put {\copy\figurelegendfive} [rt] at 3 0.8
\footnotesize
\put {$u$} [rb] at 3.0 0.04
\put {$a(u)$} [lt] at 0.025 2.2
\setsolid
\setplotsymbol ({\sixrm .})
\Black{\relax  % % square wave actual $a(u)$
\plot
    0.0022    0.3000
    0.0305    0.3001
    0.0589    0.3001
    0.0873    0.3002
    0.1157    0.3002
    0.1440    0.3001
    0.1724    0.3001
    0.2008    0.3001
    0.2292    0.3000
    0.2575    0.3000
    0.2859    0.3000
    0.3143    0.2999
    0.3427    0.2999
    0.3710    0.3000
    0.3994    0.3000
    0.4278    0.3001
    0.4562    0.3001
    0.4846    0.3001
    0.5129    0.2999
    0.5413    0.2996
    0.5697    0.2993
    0.5981    0.2990
    0.6264    0.2989
    0.6548    0.2991
    0.6832    0.2995
    0.7116    0.3004
    0.7399    0.3020
    0.7683    0.3048
    0.7967    0.3093
    0.8251    0.3168
    0.8534    0.3326
    0.8818    0.3505
    0.9102    0.3915
    0.9386    0.6602
    0.9669    1.1385
    0.9953    1.6312
    1.0237    1.9646
    1.0521    2.1312
    1.0805    2.1904
    1.1088    2.2019
    1.1372    2.2044
    1.1656    2.2033
    1.1940    2.2006
    1.2223    2.1981
    1.2507    2.1964
    1.2791    2.1954
    1.3075    2.1950
    1.3358    2.1950
    1.3642    2.1955
    1.3926    2.1962
    1.4210    2.1972
    1.4493    2.1982
    1.4777    2.1993
    1.5061    2.2002
    1.5345    2.2009
    1.5628    2.2015
    1.5912    2.2018
    1.6196    2.2018
    1.6480    2.2015
    1.6763    2.2008
    1.7047    2.1998
    1.7331    2.1984
    1.7615    2.1967
    1.7899    2.1947
    1.8182    2.1925
    1.8466    2.1903
    1.8750    2.1873
    1.9034    2.1793
    1.9317    2.1605
    1.9601    2.1249
    1.9885    2.0509
    2.0169    1.8982
    2.0452    1.6812
    2.0736    1.4589
    2.1020    1.2913
    2.1304    1.2160
    2.1587    1.1974
    2.1871    1.1980
    2.2155    1.2017
    2.2439    1.2025
    2.2722    1.2015
    2.3006    1.2000
    2.3290    1.1989
    2.3574    1.1984
    2.3858    1.1983
    2.4141    1.1986
    2.4425    1.1990
    2.4709    1.1995
    2.4993    1.2000
    2.5276    1.2003
    2.5560    1.2004
    2.5844    1.2005
    2.6128    1.2004
    2.6411    1.2003
    2.6695    1.2002
    2.6979    1.2000
    2.7263    1.1999
    2.7546    1.1998
    2.7830    1.1998
    2.8114    1.1998
    2.8398    1.1998
 /\relax}\relax

\setsolid
\setplotsymbol ({\fiverm .})
%
% plot_sq_func_iteration 1
\Blue{\relax  % iteration 1
\plot
    0.0022    0.3247
    0.0305    0.3195
    0.0589    0.3119
    0.0873    0.3006
    0.1157    0.2897
    0.1440    0.2789
    0.1724    0.2702
    0.2008    0.2625
    0.2292    0.2566
    0.2575    0.2514
    0.2859    0.2475
    0.3143    0.2444
    0.3427    0.2425
    0.3710    0.2411
    0.3994    0.2396
    0.4278    0.2386
    0.4562    0.2381
    0.4846    0.2381
    0.5129    0.2389
    0.5413    0.2402
    0.5697    0.2424
    0.5981    0.2452
    0.6264    0.2490
    0.6548    0.2533
    0.6832    0.2589
    0.7116    0.2648
    0.7399    0.2718
    0.7683    0.2787
    0.7967    0.2873
    0.8251    0.2981
    0.8534    0.3181
    0.8818    0.3539
    0.9102    0.4178
    0.9386    0.5307
    0.9669    0.7217
    0.9953    1.0103
    1.0237    1.3569
    1.0521    1.6775
    1.0805    1.8987
    1.1088    2.0029
    1.1372    2.0136
    1.1656    1.9753
    1.1940    1.9326
    1.2223    1.9141
    1.2507    1.9191
    1.2791    1.9311
    1.3075    1.9398
    1.3358    1.9445
    1.3642    1.9474
    1.3926    1.9496
    1.4210    1.9502
    1.4493    1.9483
    1.4777    1.9462
    1.5061    1.9515
    1.5345    1.9693
    1.5628    1.9952
    1.5912    2.0191
    1.6196    2.0307
    1.6480    2.0266
    1.6763    2.0115
    1.7047    1.9977
    1.7331    1.9961
    1.7615    2.0127
    1.7899    2.0460
    1.8182    2.0851
    1.8466    2.1101
    1.8750    2.1122
    1.9034    2.0939
    1.9317    2.0570
    1.9601    1.9949
    1.9885    1.9002
    2.0169    1.7778
    2.0452    1.6358
    2.0736    1.4850
    2.1020    1.3541
    2.1304    1.2581
    2.1587    1.1985
    2.1871    1.1689
    2.2155    1.1594
    2.2439    1.1599
    2.2722    1.1640
    2.3006    1.1688
    2.3290    1.1737
    2.3574    1.1781
    2.3858    1.1820
    2.4141    1.1854
    2.4425    1.1885
    2.4709    1.1910
    2.4993    1.1926
    2.5276    1.1931
    2.5560    1.1926
    2.5844    1.1914
    2.6128    1.1905
    2.6411    1.1903
    2.6695    1.1915
    2.6979    1.1940
    2.7263    1.1980
    2.7546    1.2036
    2.7830    1.2101
    2.8114    1.2172
    2.8398    1.2249
 /\relax}\relax

\Red{\relax  % iteration 4
\plot
    0.0022    0.3723
    0.0305    0.3663
    0.0589    0.3577
    0.0873    0.3447
    0.1157    0.3322
    0.1440    0.3198
    0.1724    0.3099
    0.2008    0.3011
    0.2292    0.2943
    0.2575    0.2883
    0.2859    0.2840
    0.3143    0.2804
    0.3427    0.2782
    0.3710    0.2766
    0.3994    0.2750
    0.4278    0.2739
    0.4562    0.2733
    0.4846    0.2734
    0.5129    0.2744
    0.5413    0.2760
    0.5697    0.2786
    0.5981    0.2818
    0.6264    0.2863
    0.6548    0.2914
    0.6832    0.2979
    0.7116    0.3048
    0.7399    0.3129
    0.7683    0.3210
    0.7967    0.3309
    0.8251    0.3435
    0.8534    0.3668
    0.8818    0.4084
    0.9102    0.4824
    0.9386    0.6132
    0.9669    0.8339
    0.9953    1.1672
    1.0237    1.5670
    1.0521    1.9359
    1.0805    2.1889
    1.1088    2.3062
    1.1372    2.3152
    1.1656    2.2676
    1.1940    2.2146
    1.2223    2.1891
    1.2507    2.1903
    1.2791    2.1990
    1.3075    2.2039
    1.3358    2.2040
    1.3642    2.2018
    1.3926    2.1983
    1.4210    2.1928
    1.4493    2.1839
    1.4777    2.1745
    1.5061    2.1732
    1.5345    2.1854
    1.5628    2.2064
    1.5912    2.2246
    1.6196    2.2288
    1.6480    2.2153
    1.6763    2.1897
    1.7047    2.1659
    1.7331    2.1551
    1.7615    2.1632
    1.7899    2.1887
    1.8182    2.2193
    1.8466    2.2344
    1.8750    2.2252
    1.9034    2.1946
    1.9317    2.1454
    1.9601    2.0708
    1.9885    1.9636
    2.0169    1.8300
    2.0452    1.6785
    2.0736    1.5200
    2.1020    1.3834
    2.1304    1.2835
    2.1587    1.2216
    2.1871    1.1905
    2.2155    1.1801
    2.2439    1.1799
    2.2722    1.1833
    2.3006    1.1874
    2.3290    1.1913
    2.3574    1.1948
    2.3858    1.1979
    2.4141    1.2003
    2.4425    1.2025
    2.4709    1.2042
    2.4993    1.2048
    2.5276    1.2043
    2.5560    1.2028
    2.5844    1.2007
    2.6128    1.1987
    2.6411    1.1976
    2.6695    1.1977
    2.6979    1.1992
    2.7263    1.2023
    2.7546    1.2068
    2.7830    1.2123
    2.8114    1.2183
    2.8398    1.2250
 /\relax}\relax

\endpicture
}

\begin{figure}[ht]
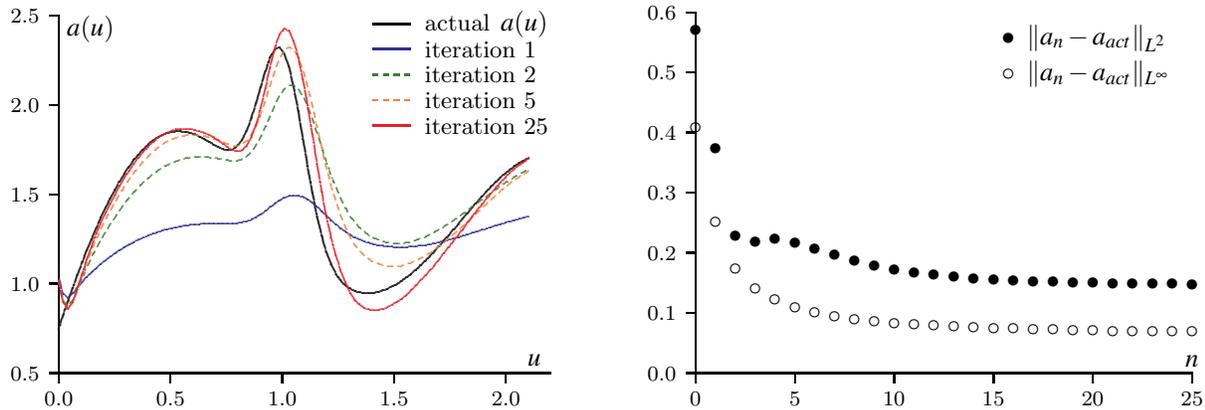

\hbox to \hsize{\hss\copy\figureone\hss\hss\copy\figuretwo\hss}
\caption{\small {\bf Recovery of $a(u)$ from time trace data}}
\label{fig:titi_titi}
\end{figure}

\begin{figure}[ht]
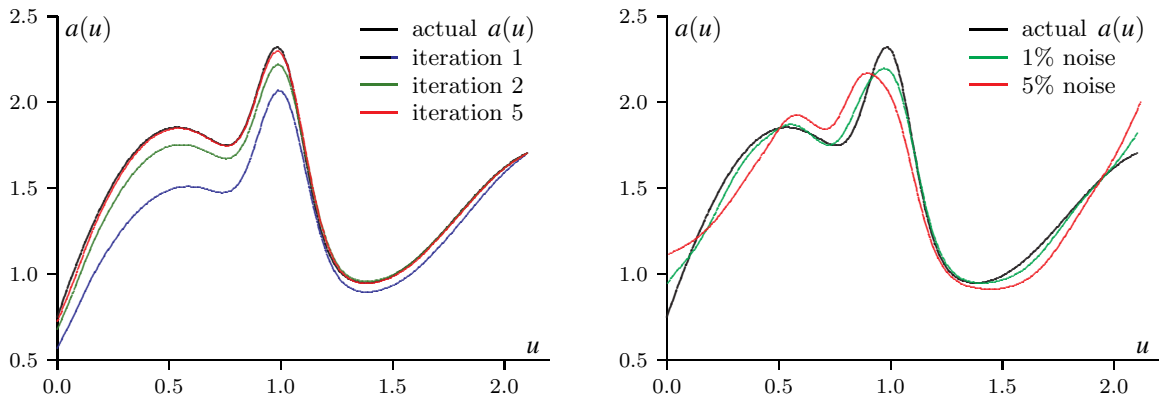

\hbox to \hsize{\hss\copy\figurethree\hss\copy\figurefour\hss}
\caption{\small {\bf Recovery of $a(u)$ from final time  data}}
\label{fig:fiti_fiti}
\end{figure}

Figure~\oldref{fig:fiti_fiti} shows reconstructions of the same 
%BK 2020-10-16 $a(u)_{act}$
$a_{act}$
but using final time $g(x) = u(x,T;a)$ measurements at $T=1$.
Again the initial approximation was $a(u) = 1$.
The leftmost figure shows the reconstructions of selected iterations
and here the fifth one 
%BK 2020-10-16 was 
corresponded to effective numerical convergence.
The noticeable difference is in both the speed and in far greater accuracy
obtained from such final time data as opposed to boundary time trace data.
This pattern was consistent in several numerical experiments with different
$a(u)$ functions and different imposed boundary and forcing functions
$r(x,t,u)$.

Note that the small level of noise added to the data $g(x)$ (at 0.1\%),
to avoid any issue of an inverse crime,
played a quite insignificant role.
The rightmost figure show the final reconstructions (again the fifth
iteration) under increased noise levels; now at 1\% and 5\%.

\begin{wrapfigure}{R}{0.5\textwidth}
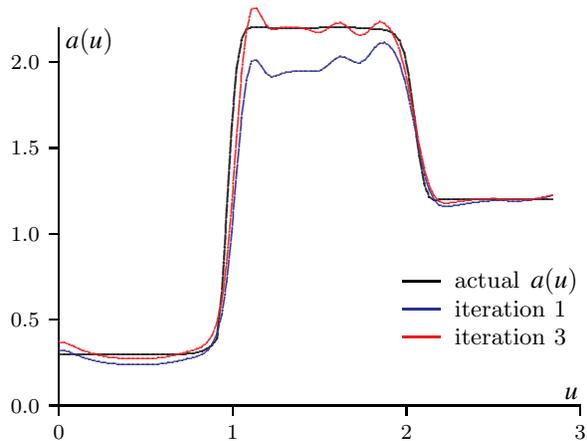

%\centering
\vspace{-18pt}
\hbox to \hsize{\hss\copy\figurefive\hss}
\caption{\small {\bf Recovery of $a(u)$ from final time  data}}
\label{fig:fiti_fiti_sq}
\vspace{-12pt}
\end{wrapfigure}

The  thermal conductivity of most materials varies considerably
over a range that encompasses a phase transition.
An experiment set up to make conductivity measurements would most likely
be conducted over a limited range sufficient to be restricted to a single
phase, but we show in Figure~\oldref{fig:fiti_fiti_sq} a simulation
of an experiment where $a(u)$ has discontinuities.

% The above \par is essential to the use of \wrapfig.
%%%%%%%%%%%%%%%%%%%%%%%%%%%%%
From an analysis viewpoint we need $a(u)$ to be in at least $C^{2+\alpha}$
and so we illustrate these boundaries with a mollified piecewise constant
function as shown.
This posed little difficulties in the case of final time data
as the figure to the right shows.
Effective numerical convergence was essentially at the third iteration
and even the first iteration starting from $a_{\rm init} = 1$ contains
much of the salient features of the actual $a(u)$.
The noise level added to the data was set at $0.1\%$ but reconstructions
of similar quality are obtained at higher noise levels than that shown in
Figure~\oldref{fig:fiti_fiti} for smoother actual $a(u)$.
The reconstruction shown is slightly under-smoothed as can be seen from
the oscillations in the region between $u\in (1,2)$.
This was done to show how accurately the algorithm was able to detect
and reconstruct sharp boundaries.

By letting our variation equal such bounds we see the efficiency of the 
numerical simulation based on the algorithm developed in
section~\oldref{subsec:fiti} and resulting in Theorem~\oldref{th:contr_fiti}.
%%%%%%%%%%%%%%%%%%%%%%%%%%%%%%%%%%%%%%%%%%%

\section*{Acknowledgment}

\noindent
The work of the first author was supported by the Austrian Science Fund {\sc fwf}
under the grants P30054 and DOC 78.

\noindent
The work of the second author was supported
in part by the
National Science Foundation through award {\sc dms}-1620138.
%%%%%%%%%%%%%%%%%%%%%%%%%%%%%%%%%%%%%%%%%%%%

\end{document}